\documentclass[a4paper,12pt]{amsart}
\usepackage{amsmath,mathrsfs,amssymb,amsthm,amscd}
\usepackage{a4wide}
\usepackage[]{fontenc}
\usepackage{hyperref}
\usepackage{graphicx}
\usepackage{color}

\usepackage{pdfsync}

\newcounter{EQNR}[NNN]
\newcommand{\nne}{\refstepcounter{EQNR}  \tag{\theNNN .\theEQNR} }
\newcommand{\nn}{\vskip 2mm 
\refstepcounter{NNN}\noindent {\bf \theNNN . }}

\DeclareMathOperator{\ca1}{\widehat{\operatorname{c}}_1}
\DeclareMathOperator{\c1}{{\operatorname{c}}_1}

\DeclareMathOperator{\dv}{div}

\DeclareMathOperator{\Nm}{Nm}

\DeclareMathOperator{\ord}{ord}
\DeclareMathOperator{\Pica}{\widehat{Pic}}

\DeclareMathOperator{\cha}{\widehat{CH}}
\DeclareMathOperator{\Spec}{Spec}

\DeclareMathOperator{\vol}{vol}

\newcommand{\fin}{\textrm{fin}}
\newcommand{\Ar}{\textrm{Ar}}

\newcommand{\dd}{\textrm{d}}

\newcommand{\can}{\textrm{can}}

\newcommand{\ZZ}{{\mathbb Z}}
\newcommand{\RR}{{\mathbb R}}

\newcommand{\CC}{{\mathbb C}}
\newcommand{\QQ}{{\mathbb Q}}

\newcommand{\HH}{{\mathbb H}}
\newcommand{\PP}{{\mathbb P}}
\newcommand{\FF}{{\mathbb F}}

\newcommand{\calD}{{\mathcal{D}}}

\newcommand{\calF}{{\mathcal{F}}}
\newcommand{\calG}{{\mathcal{G}}}

\newcommand{\calK}{{\mathcal{K}}}
\newcommand{\calL}{{\mathcal{L}}}

\newcommand{\calO}{{\mathcal{O}}}

\newcommand{\calW}{{\mathcal{W}}}
\newcommand{\calX}{{\mathcal{X}}}
\newcommand{\calY}{{\mathcal{Y}}}

\newcommand{\frakp}{{\mathfrak{p}}}

\newcommand{\belyi}{{\boldsymbol \beta}}

\newcommand{\inftyQ}{{Q}}

\def\?{\ \texttt{???}\ \immediate\write16{}%
\immediate\write16{Warning: There was still a question mark . . . }%
\immediate\write16{}}

\begin{document}

\title[Bounds for the arithmetic self-intersection 
of the dualizing sheaf]{On  the arithmetic self-intersection number 
of the dualizing sheaf on arithmetic surfaces}

\author{Ulf K\"uhn}
\email{kuehn@math.uni-hamburg.de}
\address{Department Mathematik (AZ)\\Universit\"at Hamburg\\
  Bundesstrasse 55\\D-20146 Hamburg}
\keywords{Arakelov theory, arithmetic self-intersection numbers,
  effective Mordell, modular curves, Fermat curves}    
\thanks{{\bf MSC 2000:} 14G40,11G18,11G30}

\begin{abstract} We study the arithmetic
  self-intersection number of the dualizing sheaf on arithmetic
  surfaces with respect to morphisms of a particular kind. We
  obtain  upper bounds  for the
  arithmetic self-intersection number of the dualizing sheaf on
  minimal regular models of the modular curves associated with
  congruence subgroups of the modular group $\Gamma(1)$ 
as well as the Fermat
  curves with prime exponent.
\end{abstract}

\maketitle
\tableofcontents
\section*{Introduction}
\setcounter{section}{0} \nn 
Let $X$ be a curve defined over some number field $K$ such that its 
genus
$g$ is larger than one. In the arithmetic intersection theory on a
regular model $f\!: \calX \to \Spec \calO_K$ of $X$ the hermitian line
bundle $\overline \omega_{\Ar}$, where $\overline \omega_{\Ar}$ is the
dualizing sheaf $\omega_{\calX} = \omega_{\calX/\calO_K} \otimes f^*
\omega_{\calO_K/\ZZ}$ equipped with the Arakelov metric (see
\cite{Ar}, p.1177, \cite{MB1}, p.75), plays a prominent role. Aspects
of the
arithmetic meaning of its arithmetic self-intersection number
$\overline{\omega}^2_{\Ar}$ are addressed in \cite{Sz},  \cite{Ul}
 and \cite{Zh}. 
In particular, certain conjectured upper bounds for algebraic 
families of arithmetic surfaces
(see e.g.  \cite{La}, p.166, or \cite{Sz}, p.244)
would imply an effective version of Mordell's
conjecture (cf. \cite{Pa}, \cite{MB}, or Vojta's appendix in
\cite{La}).  Thus such bounds are equivalent to the uniform
$abc$-conjecture for number fields \cite{Fr}, which in turn has many
of different applications and consequences in number theory. 
Also speculations had been made whether $\overline{\omega}^2_{\Ar}$ equals a special
 value of of the logarithmic derivative of the Hasse-Weil $L$-function of $X$.

\nn \label{assumptions} The main result of this article is a new formula for
$\overline{\omega}^2_{\Ar}$, see in particular Theorem \ref{thm:keyformula}.
As consequences of this rather technical result upper
bounds for $\overline{\omega}^2_{\Ar}$ for particular curves may be
calculated in a straightforward manner.
Prominent examples are modular curves $X(\Gamma)$
associated with (congruence) subgroups $\Gamma$ of
$\operatorname{SL}_2(\ZZ)$ and Fermat curves. Bounds for these curves
have been asked for since the beginning of Arakelov theory 
(see e.g. \cite{La}, p. 130 or \cite{MB}, 8.2). 
Except for the discrete family of modular curves $X_0(N)$  (see \cite{AU} plus \cite{MU} or \cite{JK}), there were only few results 
on estimates or even upper bounds for $\overline{\omega}^2_{\Ar}$ for other curves\footnote{Since a first version of this paper was made public H. Mayer \cite{Ma} has obtained estimates for the modular curves $X_1(N)$ by mimicking the proof of \cite{AU} plus \cite{JK}
and   more generally A. Javanpeykar \cite{Ja} was able to bound $\overline{\omega}^2_{\Ar}$ in terms of the minimal Belyi degree}.

We now set the notation for the following theorem, which will be crucial for 
the applications we have in mind. 
Let $\calY \to \Spec \calO_K$ be an arithmetic surface and 
write $Y$ for its generic fiber. Let
$\inftyQ ,P_1,...,P_r \in Y(K)$ such that $Y \setminus
\{\inftyQ,P_1, ...,P_r\}$ is hyperbolic. Then we
 consider any arithmetic surface $\calX\to \Spec \calO_K$ 
equipped with a  morphism of arithmetic surfaces 
$\belyi: \calX \to \calY$ such that the induced  morphism
$\belyi: X \to Y$  of algebraic curves
 defined over  $K$ is  unramified above
$Y(K) \setminus\{\inftyQ ,P_1,...,P_r \}$. 
 Let  $g\ge 2 $ be the genus
of $X$ and $d = \deg(\belyi)$.
We write
$\belyi^* \inftyQ = \sum b_j S_j$ and call
the points $S_j$ \emph{labeled}. Set $b_{\max} = \max_j \{b_j\}$.
Divisors on $X$ with support in the labeled points are called \emph{labeled}.
Finally, a prime $\frakp$ is said to be \emph{bad} if the fiber
of $\calX$ above $\frakp$ is reducible. Note that a prime of bad reduction need not 
be a bad prime in the above sense.

Now a consequence of our main result Theorem \ref{thm:keyformula} is the following

{\bf Theorem I.}  \label{thm:main} \emph{Let $\belyi: \calX \to \calY$ be 
a morphism of arithmetic surfaces as in \ref{assumptions}.
Assume that all labeled points  are
  $K$-rational and  that all labeled divisors of degree zero are torsion, then
  the arithmetic self-intersection number of the dualizing sheaf on  $\calX$ 
  satisfies the inequality
\begin{align}\label{eq:main}
  \overline{\omega}_\Ar^2 &\le (2 g-2) \left( \log|\Delta_{K|\QQ}|^2 + 
   [K:\QQ]  \left(  \kappa_1 \log b_{\max} + \kappa_2 \right) +
 \sum_{\frakp \,\,{\rm bad}} a_\frakp
  \log\Nm(\frakp) \right),\nne
\end{align}
where $\kappa_1= \kappa_1(Y,\inftyQ,P_1,...,P_r)\in \RR$ is a  constant
independent of $\calX$; similarly the constant  $\kappa_2$ is independent of $\calX$. 
The coefficients $a_\frakp \in \QQ$ are 
determined by certain local intersection numbers (see formula \eqref{eq:def-ap-intro} below).
}

The proof of Theorem I is an easy combination of Corollary \ref{cor:keyformula}, Theorem 
\ref{thm:keyformula}
and Proposition \ref{cor:bound-bbelyi}.

The bounds of the analytic contribution to $\overline{\omega}^2_{\Ar}$ are given by the
 constants $\kappa_1$ and $\kappa_2$, we refer to \ref{cor:bound-bbelyi} for more details.
These constants, which  come from a comparison of the hyperbolic metric on $Y$ with the Arakelov metric on $X$, are  
obtained by a refinement of results on the sup norm of
automorphic forms due to Jorgenson-Kramer \cite{JK2}. Our approach avoids
most of the technical difficulties in the previous approach by \cite{AU} plus
\cite{MU} or by \cite{AU} plus \cite{JK}. A key new ingredient
is a systematic use of the generalised arithmetic intersection theory
\cite{Kue}. Compared to the works \cite{AU} plus \cite{MU} or \cite{AU} plus
\cite{JK} this new method yields a new and considerably shorter proof
of an upper bound for the curves $X_0(N)$. Moreover it is not
restricted to the modular curves $X_0(N)$ with a square-free level $N$.

The 
geometric contribution to $\overline{\omega}^2_{\Ar}$ is reduced to a standard calculation, which is
similar to that of Abbes-Ullmo \cite{AU}. In addition we were able to give upper bounds, 
which can be computed quite easily whenever the model is known.
To keep the notation simple, we write $S_j$ also for the Zariski
closure in $\calX$ of a labeled point $S_j$.  Let $\calK$ be a canonical
divisor of $\calX$, then for each labeled point $S_j$ we can find a divisor
$\calF_j$ such that
\begin{align*}
 \left( S_j + \calF_j -\frac{1}{2g-2} \calK\right ). \,\mathcal{C}^{(\mathfrak{p})}_{l}=0
\end{align*}
for all irreducible components $\mathcal{C}^{(\mathfrak{p})}_{l}$ of
the fiber $f^{-1}(\mathfrak{p})$ above $\mathfrak{p} \in \Spec
\mathcal{O}_K$. 
Similarily we find  
 for each labeled point $S_j$   a
divisor $\calG_j$ such that also for all $ \mathcal{C}^{(\mathfrak{p})}_{l}$ as before
\begin{align*}
\left ( S_j + \calG_j -\frac{1}{d} \belyi^*\inftyQ \right) .\, \mathcal{C}^{(\mathfrak{p})}_{l}=0.
\end{align*}

Then the rational numbers $a_\frakp$ in the theorem are determined by 
the following arithmetic intersection numbers
of trivially metrised hermitian line bundles 
\begin{align}\label{eq:def-ap-intro}
 \sum_{\frakp \,\,{\rm bad}} a_\frakp
  \log\Nm(\frakp)
&= - \frac{2 g 
}{d} \sum_j b_j\,\calO( \calG_j)^2 +
\frac{
2g-2}{d} \sum_j b_j\, \calO(\calF_j)^2.  \nne 
\end{align}

\nn Since the calculations of the quantities $a_\frakp$ in Theorem I
are rather involved, we consider instead for each bad prime $\frakp$ the more easy accessible quantity
$b_\frakp$, which is defined as follows. Let
\begin{align*}
\calX \times \overline
\FF_\frakp
= \sum_{j=1}^{r_\frakp} m_j C_j^{(\frakp)}
\end{align*}
be the decomposition into irreducible components
and set 
\begin{align*}
u_\frakp= \max_{i,j}
|C_i^{(\frakp)}.C_j^{(\frakp)}|, \qquad
l_\frakp= \min_{C_i^{(\frakp)}.C_j^{(\frakp)}\neq 0}
|C_i^{(\frakp)}.C_j^{(\frakp)}|.
\end{align*}
We further denote by $c_\frakp$ the 
connectivity of the dual graph of $\calX \times \overline
\FF_\frakp$, i.e. the 
minimal number
of intersection points needed to connect any two
irreducible components. 
Using these information we can define
\begin{align*}
 b_\frakp = \left( 
 \sum_{k=1}^{c_\frakp} 
\left( \sum_{l=1}^{ k} 
 \left(\frac{u_\frakp}{l_\frakp}\right)^{l-1}  \right)^2
+ (r_\frakp - c_\frakp-1) \left( \sum_{l=1}^{ c_\frakp} 
 \left(\frac{u_\frakp}{l_\frakp}\right)^{l-1}  \right)^2 \right)  
 \frac{u_\frakp}{l_\frakp^2} 
\end{align*}

{\bf Theorem II.} \emph{The rational numbers  $a_\frakp$ in Theorem I  satisfy 
$a_\frakp \le 2 g\, b_\frakp$.}

\nn {\bf Applications.} If $\calX$ is a suitable model of a modular curve 
associated with a congruence subgroup $\Gamma$,
then  because of the Manin-Drinfeld theorem (see e.g. \cite{El}) the 
assumptions of Theorem I are satisfied  for the natural morphism given by 
the $j$-map. 
In particular, if $\Gamma$ is of a certain kind, then the coefficients
 $a_\frakp$ in \eqref{eq:main} can be computed explicitly by means of
 the descriptions of models for $X(\Gamma)$ (see e.g. \cite{KM},
 \cite{DR}). We illustrate this with the following theorem.

{\bf Theorem III.}
\emph{Let
$N$ be a square free integer
   having at least two different prime factors and $(N,6)=1$.
Let $\calX_0(N)$ be the minimal regular model
   of the modular curve $X_0(N)$ and $g_N$ its genus. Then the
   arithmetic self-intersection number of its dualizing sheaf equipped
   with the Arakelov metric is bounded from above by
\begin{align*}
\overline{\omega}_\Ar^2  
&\le    (16 \pi \kappa_{\circ}-1) \, g_N \log(N) +f(N) 
\end{align*}
where $\kappa_{\circ} \in \RR$ is the Jorgenson-Kramer constant for $\operatorname{SL}_2(\ZZ)$ and
$f(N)$ is an explicit function of order $O(g_N)$.
}

We point to the fact that
our upper bound  together with a relatively simple lower bound as in \cite{AU} 
imply the asymptotic $\overline{\omega}_\Ar^2 \sim g_N\log(N)$. 
In contrast to the previously known assymptotic  
\begin{align}\label{thm:au-mu}
\overline{\omega}_{\calX_0(N),\Ar}^2= 3 g_N \log(N) (1+O(\log\log(N)/
\log(N)))\nne
\end{align}
this slightly weaker assymptotic 
is easily deduced from the results in this paper.
With methods, which depend strongly on the specific arithmetic of
$\Gamma_0(N)$,   Abbes and Ullmo  derived in \cite{AU}
  two different formulae for
$\overline{\omega}_{\calX_0(N),\Ar}^2$, which contained certain
quantities that they  were unable to estimate from above. In \cite{MU} Michel and Ullmo
provided estimates for certain integrals of Eisenstein series against
the canonical volume form. These estimates together with a calculation
of the Neron-Tate height of certain Heegner divisors proved the formula
\eqref{thm:au-mu} by using one of the formulae of \cite{AU}.  In
\cite{JK} Jorgenson and Kramer obtained estimates for derivatives of
the Selberg zeta functions, which  also yields this estimate by means of the other formula of
\cite{AU}. However both of this two approaches need strong arithmetic facts on $\Gamma_0(N)$, 
which are not needed for the slightly weaker assymptotic.

For the next results the 
quantities $b_\frakp$ are used,

{\bf Theorem IV.} \emph{Let $N$ be an integer having at least two
  different prime factors and $(N,6)=1$.  Let $\calX(N)$ be the regular model 
of the modular curve $X(N)$ given  by the moduli description (see \cite{KM}).  Then the
  arithmetic self-intersection number of its dualizing sheaf equipped
  with the Arakelov metric is bounded from above by
\begin{align*}
  \overline{\omega}_\Ar^2 &\le (2 g_N-2) \left( \log|\Delta_{\QQ(\zeta_N)|\QQ}|^2 + 
   [\QQ(\zeta_N):\QQ]  (  \kappa_1 \log N +\kappa_2)   
 + 2  \sum_{\frakp \supset (N)} b_\frakp \log \Nm \frakp  
  \right),
\end{align*} 
where  $g_N$ denotes the genus of $X(N)$, 
$\kappa_1, \kappa_2 \in \RR$ are constants independent of $N$ and 
\begin{align*} 
 b_\frakp \le \frac{(r_\frakp-1)^2 m_\frakp}{s_\frakp}, 
\end{align*}
where $r_\frakp$ is as above, $m_\frakp \le p^{2 v_{p}(N)-2}$ where $(p)=\frakp \cap\ZZ$ and   
$s_\frakp$ is the number of supersingular points above $\frakp$.}

Other examples of curves where our result could be applied are the
Fermat curves. Here we consider just the Fermat curves
with prime exponents.

{\bf Theorem V.}
\emph{Let $\calX$ be the desingularisation of the
  closure in $\PP^2_{\ZZ[\zeta_p]}$ of the Fermat curve $x^p+y^p=z^p$ 
with prime exponent $p$ (see \cite{Mc}).  
Then the arithmetic self-intersection
  number of its dualizing sheaf equipped with the Arakelov metric is
  bounded from above by
\begin{align*}
  \overline{\omega}_\Ar^2 &\le (2 g_p-2) \left( \log|\Delta_{\QQ(\zeta_p)|\QQ}|^2 + 
   [\QQ(\zeta_p):\QQ]  (  \kappa_1 \log p +\kappa_2)   
 +  g_p\, p^8 \log p     
  \right),
\end{align*}
where $\kappa_1, \kappa_2 \in \RR$ are  constants independent of $p$.}

In  forthcoming papers
\cite{C},\cite{CK} we explicitly calculate the coefficients $a_\frakp$ 
for Fermat curves and other modular curves.

We finally observe that there are only finitely many isomorphism
classes of curves provided with a morphism of a fixed degree and fixed
branch points.  Therefore no continuous family of curves satisfying
the assumptions of Theorem I exists. However, one may hope that 
the ideas of this article may be helpful for exact calculations of $\overline{\omega}_\Ar^2$, e.g.  
by extending Zagiers calculation in \cite{Za}. For the moment the speculation whether
$\overline{\omega}_\Ar^2$ equals a special value of the logarithmic derivative of the Hasse-Weil $L$-function has to be considered completly open, cf. \cite{Zh2}.

\nn {\bf Plan of paper.}  In the first and second sections of this
paper we present the necessary background material on Arakelov theory
on arithmetic surfaces.  After this preparatory work we study in the
third section the behavior of arithmetic intersection numbers with
respect to a finite morphism.  The fourth section is devoted to the
analytical aspects needed in our bound \eqref{eq:main}. In section
five we study the geometric aspects used in the bound of the
quantities $a_\frakp$ of \eqref{eq:main}.  In the final three sections
we apply our result to the above mentioned examples.

 \section{Basic properties of Green functions.}\label{sec:basic-green} 
 We recall here some basic facts related to normalised Green
 functions on a compact Riemann surface $X$. This is mainly due to the
 fact that inconsistent normalizations for these Green functions are
 used in the literature.  Our normalization is as in \cite{SABK},
 \cite{La} and \cite{Kue}; it is twice the Green function advocated
 in \cite{BKK1} and it is $-2$ times the normalization used in
 \cite{El2}, \cite{Sz} and \cite{Fa}.

 \nn {\bf Volume forms.}  Let $\nu$ be a volume form (cf. \cite{La}
 p.21), i.e. a smooth, positive, real $(1,1)$-form with $\int_X
 \nu(z)=1$. Given two volume forms $\mu, \nu$ we can write $\mu(z) =
 f(z) \nu(z)$, where $f$ is a smooth function. 

 If $\belyi : X \to Y$ is
 a morphism of compact Riemann surfaces and if $\mu$ is a volume form
 on $Y$, then $(\deg \belyi)^{-1} \,\belyi^* \mu$ need not be
 volume form anymore. Indeed, the ramification forces the pull-back
 form to vanish. On the other hand, it may happen that the pull-back of 
a singular form becomes a volume form.

If the genus of  of a Riemann surface $X$ is greater than one, then  we
have   the \emph{canonical volume form}
\begin{align*}
\nu_\can (z) = \frac{i}{2g} \sum_j |f_j|^2 \dd z \land \dd \overline z,
\end{align*}
where $f_1(z) \dd z$, ..., $f_g(z) \dd z$ is an
orthonormal basis of $H^0(X,\Omega^1)$ equipped with the
natural scalar product.

 \nn {\bf Green functions.} \label{def:greenfunc} A \emph{Green function} $g$
 associated with $\nu$ is a real valued function on $X\times X$ which
 is smooth outside the diagonal and has    an expansion
 $g(z,w) = - \log|z-w|^2 + h(z,w)$ with a smooth function $h$ near the diagonal. As a
 current it satisfies
\begin{align}\label{eq:greeneq}
\operatorname{dd}^c\, [g(w,z)] + \delta_w =[\nu(z)]\nne
\end{align}
and it is called \emph{normalized} if in addition for all $w \in X$
\begin{align}\label{eq:greennormal}
\int_X g(z,w) \nu(z) = 0.\nne
 \end{align}
One can show that there is a unique normalized Green function $g_\nu$
associated to $\nu$.

\nn {\bf Resolvent Kernel.} Normalized Green functions are also
referred to as the resolvent kernel. In order to describe this, we
consider the space $\mathcal{C}^\infty(X,\CC) \subseteq  L^2(X,\nu)$ of complex valued
$\mathcal{C}^\infty$-functions on $X$.  Then the Laplace operator
$\Delta=\Delta_\nu$ associated with a volume form $\nu$ is defined by
\begin{align} \label{eq:laplaceddc}
\Delta f \cdot \nu = \frac{1}{\pi i} \partial \bar \partial f = -
\operatorname{dd}^c f,\nne
\end{align}
here $f \in \mathcal{C}^\infty(X,\CC)$.
The eigenvalues of $\Delta$ are positive real numbers  
$
0 = \lambda_0 < \lambda_1 \le \lambda_2 \le \ldots \le \lambda_m \le
\ldots
$ 
satisfying $\lim_{m \to \infty} \lambda_m = \infty$.
We let $\phi_0=1$, $\phi_1$, $\ldots$ denote the corresponding 
normalized eigenfunctions, i.e., 
$
\int_X \phi_n \overline \phi_m \nu = \delta_{n,m}
$.
If  $f \in \mathcal{C}^\infty(X,\CC)$, then
$f$ has an expansion $f(z)= \sum_{m\ge 0} a_m \phi_m(z)$, where $a_m =
\int_X f(u) \overline{\phi}_m(u) \nu(u)$. With the above normalisation 
we have the following results. 

\nn {\bf Theorem.} \label{thm:resolvent} \emph{For $f(z)= \sum_{m\ge
    0} a_m \phi_m(z) \in \mathcal{C}^\infty(X,\CC)$ we
  set 
\begin{align*}
\tilde G_\nu(f) (z)= \sum_{m>0} \frac{a_m }{\lambda_m} 
\phi_m(z)
. \end{align*}
Then we have
\begin{align*}
\tilde G_\nu(f)(z)  = \int_X g_\nu(z,w) f(w) \nu(w). 
\end{align*}
}

{\bf Proof.} This fact is well-known, 
see e.g. \cite{El2}, p.94, or \cite{Fa}, p.394. For the convenience of
the reader we recall the proof for our normalisation.
 It suffices to show that for all $f\in 
 \mathcal{C}^\infty(X,\CC)$ such that $\int_X f(z) \nu(z)=0$ we have the
 equality
 \begin{align*}
  f(P) = \int_X g_\nu(P,z) \Delta f(z) \nu(z), 
 \end{align*}
 since the kernel of $\Delta$ is spanned by the constant function $1$.
 Using the identity \eqref{eq:laplaceddc} and
  the Green's equation
 \eqref{eq:greeneq} we derive
 \begin{align*}
  \int_X g_\nu(P,z) \Delta f(z) \nu(z)&=
 -  \int_X g_\nu(P,z)\operatorname{dd}^c f(z) =
 -\operatorname{dd}^c[g_\nu(P,z)](f)  \\
 &=  \delta_P(f) - [\nu](f) =  f(P).
 \end{align*}
  \hfill $\Box$

\nn {\bf Corollary.}\label{cor:resolvent} \emph{If in addition $\int_X
  f(z) \nu(z) =0$, then $\operatorname{dd}^c -\int_X g_\nu(z,w) f(w)
  \nu (w)= f(z)\nu (z)$.}

{\bf Proof.} 
 By Theorem \ref{thm:resolvent} we have 
$\Delta \tilde G_\nu(f)(z) =  f(z)$. Using the identity
\eqref{eq:laplaceddc}
we deduce 
\begin{align*}
\operatorname{dd}^c \left( -\int_X g_\nu(z,w) f(w) \nu (w) \right)
=  \Delta \tilde G_\nu(f)(z) \cdot \nu(z) = f(z) \nu(z).
\end{align*}
\hfill $\Box$

\nn {\bf Corollary.} \label{cor:gaac}\emph{Let $g_\nu(z,w)$, $g_\mu(z,w)$ 
be the normalized Green functions associated with the volume forms 
$\nu$, $\mu$, respectively. Then
\begin{align}
 g_\nu(z,w) = g_\mu(z,w) + a_{\mu,\nu}(z) + a_{\mu,\nu}(w) + c_{\mu,\nu},\nne
\end{align}
where for $z \in X$ we set
\begin{align}
a_{\mu,\nu}(z) &= -\int_X g_\mu(z,u) \nu(u),\nne \\
c_{\mu,\nu} &= -\int_X
a_{\mu,\nu}(w) \nu(w) = \int_{X\times X}   g_\mu(w,u) \nu(u) \nu(w).\nne
\end{align}
}

{\bf Proof.} The right hand side is smooth except for a logarithmic singularity
along the diagonal in $X\times X$. It is also orthogonal to
$\nu$ in both variables. The positivity assumption on $\mu$ allows us to write
$\nu(u)= (f(u) +1)\cdot \mu(u)$ with $f \in 
\mathcal{C}^\infty(X,\CC)$. Observe that $\int_X f(u) \mu(u) =0$ and that
since $g_\mu$ is normalized, we have $a_{\mu,\nu}(z)= -\int_X
g_\mu(z,u) f(u) \mu(u)$. From corollary \ref{cor:resolvent} we deduce
$\operatorname{dd}^c a_{\mu,\nu}(z) = f(z) \mu(z) =\nu(z) -\mu(z)$,
hence the right hand side equals $g_\nu$ by uniqueness.  \hfill $\Box$

\nn {\bf Lemma.} \label{lem:intsecond} 
\emph{Let $\nu$, $\mu$  be volume forms on $X$. Then for any $P \in X$ we have
\begin{align*}
  \int_{X}
 \left( g_\nu(z,P) - g_\mu(z,P)\right)\cdot \left( \mu(z) + \nu(z) \right) = 2
 \, a_{\mu,\nu}(P) + c_{\mu,\nu}.
\end{align*}}
 
{\bf Proof.} By means of corollary \ref{cor:gaac} we have
$
   g_\nu(z,P)-g_\mu(z,P)
= a_{\mu,\nu}(z) +a_{\mu,\nu}(P) + c_{\mu,\nu} 
$.
By changing the order of integration we derive the identity
\begin{align*}
\int a_{\mu,\nu}(z) \left( \mu(z) + \nu(z) \right)
&= \int \left( \int - g_\mu (z,w) \nu(w) \right)\left( \mu(z) + \nu(z) \right) 
= -c_{\mu,\nu},
\end{align*}
thus, since $\int_X \left( \mu(z) + \nu(z) \right)=2$, 
 we obtain the claim. \hfill $\Box$ 

 \nn {\bf Lemma.} \label{lem:lowbound} \emph{Let $\nu$, $\mu$ be
   volume forms on $X$ and write $\nu(z)= (f(z)
  +1)
\cdot \mu(z)$. Then we have
\begin{align}\label{eq:bound-c-smooth}
0 \le  c_{\mu,\nu}\le \frac{2}{\lambda_1} ||f||^2_{L^2},\nne
\end{align}
where $\lambda_1$ is the first non zero eigenvalue of $\Delta_\mu$ and
$ ||f||^2_{L^2}= \int_X f(z) \overline{f(z)} \mu(z)$.
}

{\bf Proof.} We first observe that $f(z)$ is a smooth, square integrable function that is orthogonal 
to the constants with respect to $\mu$. By Theorem \ref{thm:resolvent} we have, since 
$f(z)=\overline{f(z)}$, the upper bound
\begin{align*} 
  c_{\mu,\nu} = \int_X \tilde G_\mu(f)(z) \overline{f(z)} \mu(z)
  = \sum_{m>0} 2 \frac{|a_m|^2}{\lambda_m} \le \frac{2}{\lambda_1}
  \sum_{m>0} |a_m|^2 = \frac{2}{\lambda_1}\int_X f(z) \overline{ f(z)}
  \mu(z),
\end{align*}
where as above $f(z)= \sum_{m\ge
    0} a_m \phi_m(z)$.
The lower bound is  obvious.
 \hfill $\Box$

 \section{Hyperbolic Green function}

In the sequel we need generalisations of the
 basic properties of Green functions associated with smooth volume forms
presented in section \ref{sec:basic-green}. We have to consider Green functions
associated with the hyperbolic metric as well. This metric is logarithmically singular in the sense of  \cite{Kue}.

\nn {\bf Hyperbolic curves over $\CC$.} \label{CCModelle}
Let $\HH:=\{ z = x+ i y \in \CC \, | \,  {\rm Im}\, z = y > 0\}$ denote 
the  upper half plane.
We fix a fuchsian
 subgroup $\Gamma_{\!\!\circ} \subset {\rm PSL}_2(\RR)$ of the first kind. Then 
for any subgroup $\Gamma$  of finite index  in $\Gamma_{\!\!\circ}$ the
quotient $\Gamma \setminus \HH$ by the natural action of $\Gamma$ on
$\HH$ has the structure of a Riemann surface, which
can be compactified by adding finitely many cusps. We denote this compact 
Riemann surface by $X(\Gamma)$ and
we call $X(\Gamma)$ a \emph{hyperbolic curve}.
 
On $X(\Gamma)$ there may  also be finitely many elliptic points; those
correspond to $z \in \HH$ which have a non-trivial isotropy subgroup in
$\Gamma$.  It is common to identify a local coordinate $z$ on $X(\Gamma)$ 
with a preimage $z \in \HH$.  

  The inclusion of a subgroup $\Gamma$ of
finite index in $\Gamma_{\!\!\circ}$ induces a finite, holomorphic map $ f:
X(\Gamma) \longrightarrow X(\Gamma_{\!\!\circ}) $ of degree ${\rm deg}\, f =
[\Gamma_{\!\!\circ}:\Gamma]$.  This morphism of compact Riemann surfaces is
ramified only above the  elliptic points and the
cusps of $X(\Gamma_{\!\!\circ})$.  For more details, cf. \cite{Sh}.

\nn {\bf Hyperbolic volume form on $X(\Gamma)$.} The $(1,1)$-form 
associated to the hyperbolic metric on $\HH$ is given by 
\begin{align}\label{eq:defmuhyp}
\mu= \frac{dx\wedge dy}{y^2}= \frac{i}{2} \cdot 
\frac{dz \wedge d\bar z}{ ({\rm Im} \,z)^2}.\nne
\end{align}
By abuse of notation we denote the $(1,1)$-form induced by $\mu$ on
$X(\Gamma)$ also by $\mu$. It is singular and its singularities occur
at the elliptic fixed points and at the cusps of $X(\Gamma)$ (see
e.g. \cite{Kue}, p.222).
We write $\mu_{\Gamma}$
for the log-log singular volume form on $X(\Gamma)$ determined by
$\mu$;
i.e.  $\int_{X(\Gamma)} \mu_\Gamma=1$. We call $\mu_\Gamma$ a \emph{hyperbolic volume form}.

The following facts
are crucial for our applications.
If $\nu$ is a volume form on $X(\Gamma)$, then we can write 
$\nu = f(z) \mu_\Gamma$, where $f(z)$ is a smooth function on $X(\Gamma)$. 
If $f: X(\Gamma) \to X(\Gamma_{\!\!\circ})$ is a morphism of hyperbolic curves 
induced by $\Gamma \subset \Gamma_{\!\!\circ}$, then
$f^* \mu_{\Gamma_{\!\!\circ}} = \deg(f) \mu_\Gamma$.

\nn {\bf Sup norm bounds.} Let $\Gamma \subseteq \Gamma_{\!\!\circ}$
be a  subgroup of finite index $d=[\Gamma_{\!\!\circ}:\Gamma]$
and suppose that the genus $g$ of $X(\Gamma)$ is at least $2$.  It is shown by
Jorgenson and Kramer in \cite{JK2}, Theorem 3.1, that there exists a
constant $\kappa_\circ$ only depending on $X(\Gamma_{\!\!\circ})$ such
that everywhere on  $X(\Gamma)$ 
  the estimate
\begin{align}
\label{eq:thmJoKr}
\nu_\can \le \kappa_\circ \vol(X(\Gamma_{\!\!\circ})) \frac{d}{g} 
\mu_{\Gamma}
\nne
\end{align}
holds, where $\vol(X(\Gamma_{\!\!\circ}))$ equals the volume of $X(\Gamma_{\!\!\circ})$
 with respect to the hyperbolic measure \eqref{eq:defmuhyp}.
In Theorem \ref{thm:jk+} below we consider a more precise  
comparison of these forms  that allows us to bound a certain integral.

\nn {\bf Spectral theory.} The space of  square integrable functions with
respect to $\mu_\Gamma$ on $X(\Gamma)$ will be denoted by
$L^2(X(\Gamma),\mu_\Gamma)$. Any smooth differentiable function $f(z)
\in L^2(X(\Gamma),\mu_\Gamma)$ has the following spectral
decomposition
\begin{align*}
f(z)= \sum_{n\ge 0} a_n \phi_n(z) + \sum_{\kappa} \int_0^\infty h_\kappa(t)
E_\kappa(z, \frac{1}{2}+i t) {\rm d}t, 
\end{align*}
with $\{\phi_n(z)\}$ an orthonormal basis of eigenfunctions for the
discrete
spectrum of the hyperbolic laplacian $\Delta_\Gamma= y^2
\left(\frac{\partial^2}{\partial x^2} + \frac{\partial^2}{\partial
    y^2} \right)$
and $\{ E_\kappa\}$ a complete set of Eisenstein series for the cusps
$\kappa$ of $\Gamma$. One has 
\begin{align*}
\|f(z)\|_{L^2}^2= \sum_{n\ge 0} |a_n|^2  + 
2 \pi \sum_{\kappa} \int_0^\infty |h_\kappa(t)|^2 {\rm d}t. 
\end{align*} 

\nn {\bf Hyperbolic Green function.}
The \emph{automorphic Green function} is defined for
$\text{Re}\,s >1$ by the convergent series
$$
g_\Gamma(z_1,z_2;s) := 2 \pi \cdot \sum_{\gamma \in \Gamma} - 2
Q_{s-1} \left( 1 + \frac{ |z_1-\gamma z_2|^2}{2 y_1 {\rm Im}(\gamma z_2)}
\right) .
$$
A well-known result, see e.g.\cite{Fay} Thm. 2.3, states that we
have meromorphic continuation in $s$ and that in $s=1$  the
Laurent expansion is given by
$$
g_\Gamma(z_1,z_2,s) = \frac{3/\pi}{[\Gamma:\Gamma(1)]} \frac{1}{
  s(s-1)} + g_\Gamma(z_1,z_2) + O(s-1),
$$
where the $O(s-1)$ term is smooth in $z_1,z_2$.  We call the
constant term at $s=1$ the \emph{hyperbolic Green function} for
$X(\Gamma)$.  By construction, the function $g_\Gamma(z_1,z_2)$ is
$\Gamma$-invariant and therefore descends to a function on
$X(\Gamma)$.  It is well-known that $g_\Gamma(z_1,z_2)$ has a
logarithmic singularity along the diagonal and that it is smooth
outside the cusps and the elliptic fixed points, where it has some mild
loglog singularities.  Another construction of it by means of
differentials of the third kind was given by W. Roelke \cite{Roe}, p.~22. He
showed that $g_\Gamma(z_1, z_2)$ is the normalised
Green function for $\mu_\Gamma$ whenever $z_1$ and $z_2$ are not cusps of $X(\Gamma)$.

\nn {\bf Lemma.} \label{lem:gaac-hyp} \emph{Let $g_\nu(z,w)$
be the normalized Green function associated with a volume form 
$\nu$ and let $g_\Gamma(z,w)$ be the hyperbolic Green function. Then
if $P$ is a cusp of $X(\Gamma)$ we have    the equality
 \begin{align*}
 g_\nu(z,P) = \widetilde g_\Gamma(z,P; w_0) + a_{\mu_\Gamma,\nu}(z) + 
\widetilde a_{\mu_\Gamma,\nu; w_0}(P) + c_{\mu_\Gamma,\nu}
\end{align*}
for the Green function 
considered as a function in $z$.
Here we  choose a point $w_0$ away from the cusps and set
\begin{align*}
\widetilde g_\Gamma(z,P; w_0) &= \lim_{ w \to P} \left(g_\Gamma(z,w)
  -g_\Gamma(w_0, w) \right),\\
   a_{\mu_\Gamma,\nu}(z) &= -\int_X g_\Gamma(z,u) \nu(u), \\
\widetilde a_{\mu_\Gamma,\nu, w_0}(P) &= \lim_{ w \to P}
\left(a_{\mu_\Gamma,\nu}(w) + g_\Gamma(w_0, w) \right),\\ 
c_{\mu_\Gamma,\nu} &=  \int_{X\times X}   g_\Gamma(v,w) \nu(w) \nu(v).
\end{align*}
}

{\bf Proof.} It is well known that corollary \ref{cor:resolvent} also holds
for $\mu_\Gamma$ instead of a volume form $\nu$, see e.g. \cite{Roe}, p.~35. 
From this it is clear that analogous to  \ref{cor:gaac} we have, 
with notation as above,  
\begin{align*}
 g_\nu(z,w) = g_\Gamma(z,w) + a_{\mu_\Gamma,\nu}(z) + 
a_{\mu_\Gamma,\nu}(w) + c_{\mu_\Gamma,\nu}
\end{align*}
for $z, w \in X(\Gamma) \setminus \{ {\rm cusps}\}$ .
By uniqueness of $ g_\nu(z,w)$, the right hand side is a smooth function
for all pairs $z$, $w$ outside the diagonal.  The loglog singularity
of $g_\Gamma(z,w)$ for $w$ approaching a cusp $P$ is independent of
$z$, see e.g. \cite{Fay}, p.161. Therefore the functions $\widetilde g_\Gamma(z,P; w_0) $ and
$\widetilde a_{\mu_\Gamma,\nu, w_0}(P)$ are well-defined.  Again   the claim 
follows by the
uniqueness of $ g_\nu(z,w)$.  \hfill $\Box$

\nn {\bf Remark.} 
In \cite{AU} versions of Theorem \ref{thm:resolvent}, its Corollaries
 \ref{cor:resolvent}, \ref{cor:gaac} and Lemma \ref{lem:lowbound} are
 used to express the Green function for the canonical metric evaluated
 at two different cusps in terms of the Green function for the hyperbolic
 metric on the modular curve $X_0(N)$.  In this article we only have 
  to consider the evaluation at one
 cusp instead of two, thus simplifying the arguments in \cite{AU}.

\nn {\bf Lemma.} \label{lem:intsecond-hyp} \emph{Let $\nu$ be a
  volume form on $X(\Gamma)$, then for a cusp $P \in X(\Gamma)$ we have
\begin{align*}
  \int_{X} \left( g_\nu(z,P) - \widetilde g_\Gamma(z,P;
    w_0)\right)\cdot \left( \mu_\Gamma(z) + \nu(z) \right) = 2 \,
  \widetilde a_{\mu_\Gamma,\nu; w_0}(P) + c_{\mu_\Gamma,\nu}.
\end{align*}} 

{\bf Proof.} By means of Lemma \ref{lem:gaac-hyp} the claim follows 
along the lines of the proof of Lemma \ref{lem:intsecond}.
\hfill $\Box$

 \nn {\bf Lemma.} \label{lem:bound-c-hyp} \emph{Let $\nu$  be a
   volume form on $X(\Gamma)$ and write $\nu(z)= (f(z)
  +1)
\cdot \mu_\Gamma(z)$, then we have
\begin{align}\label{eq:bound-c-hyp}
0 \le  c_{\mu_\Gamma,\nu}\le \frac{2}{\min(\lambda_1,1/4) } ||f||^2_{L^2},\nne
\end{align}
where $\lambda_1$ is the first non zero eigenvalue of $\Delta_{\Gamma}$ and
$ ||f||^2_{L^2}= \int_X f(z) \overline{f(z)} \mu_\Gamma(z)$.
}

{\bf Proof.} We first observe that $f(z)$ is a smooth, square integrable
function that is orthogonal to the constants with respect to
$\mu$. Now with the same arguments as in Lemma \ref{lem:lowbound} (see
also \cite{AU}, p.66) we derive the the upper bound
\begin{align*} 
  c_{\mu_\Gamma,\nu} &=  \int_{X(\Gamma)\times X(\Gamma)} 
 g_\Gamma(z,w) f(z) \overline{f(w)}\mu_\Gamma(w) \mu_\Gamma(z)\\
 & =2 \left(  \sum_{m>0}  \frac{|a_m|^2}{\lambda_m} + 2 \pi \sum_{\kappa} 
 \int_0^\infty\frac{| h_\kappa(t)|^2}{\frac{1}{4} + t^2 } \operatorname{d}t\right)\\
 &\le 2 \left( \frac{1}{\lambda_1}
  \sum_{m>0} |a_m|^2  +  \frac{2\pi}{\frac{1}{4}} \sum_{\kappa} 
 \int_0^\infty |h_\kappa(t)|^2 \operatorname{d}t \right)\le
 \frac{2}{\min(\lambda_1,\frac{1}{4}) } ||f||^2_{L^2}
\end{align*}
and again the lower bound $0 \le c_{\mu_\Gamma,\nu}$ is obvious.
 \hfill $\Box$

\nn{\bf Theorem.}\label{thm:jk+} \emph{Let $\Gamma \subseteq \Gamma_{\!\!\circ}$ be a
  subgroup of finite index and suppose that the genus $g$ of $X(\Gamma)$ is at least $2$. Let $f : X(\Gamma) \to
  X(\Gamma_{\!\!\circ})$ be the natural morphism of degree $d=[\Gamma:
  \Gamma_{\!\!\circ}]$.  
Let $L$ be a line bundle on $ X(\Gamma_{\!\!\circ})$
equipped with either a smooth hermitian or a $\mu_{\Gamma_{\!\!\circ}}$-admissible 
metric\footnote{cf. Definition \ref{def:hypadmissible} }.
Let $l$ be a section of  $L$ 
   whose divisor $\dv(l)$ is effective. Then there exists a constant
  $\kappa \in \RR$ independent of $\Gamma$, which depends only on
  $\log\|l\|$ as a function on $X(\Gamma_{\!\!\circ})$, such that
\begin{align} \label{bound:smoothint2}
\int_{X(\Gamma)} -\log\|l\| \, \nu_\can  \le 
 4 \pi \kappa_\circ \sum_{S_j {\rm\, cusp} }  \ord_{f(S_j)}(l) \frac{b_j \log(b_j)}{g} +  \frac{d}{g} \kappa
,\nne
\end{align}
where $\nu_\can$ is the 
canonical volume form on $X(\Gamma)$, $b_j$ denotes the ramification index 
at the cusp $S_j$ and $g$ is the genus of
$X(\Gamma)$. 
}

{\bf Proof.}  
We first assume that the divisor of $l$ is disjoint from the cusps. Then
$\|l\|$ is either a smooth function or a logarithmically-singular
function.  Thus we may apply inequality \eqref{eq:thmJoKr} to the
integral under consideration and we easily derive the upper bound.

However, if
$l$ vanishes at a cusp of $X(\Gamma)$, then the logarithmic
singularities of $\log\|\l\|^2$ together with the singularities of
$\mu_\Gamma$ near the cusps prevent us from applying inequality
\eqref{eq:thmJoKr} directly. The singularities at the elliptic points
do not matter.

Without loss of generality we may assume that $X(\Gamma_{\!\!\circ})$
has only one cusp, which we denote by $\inftyQ$; the general case
of several cusps can be proved using analogous arguments. 

For an small $\varepsilon$ we set
$$
X(\Gamma)_{\varepsilon}  = X(\Gamma) 
\setminus f^* B_\varepsilon(\inftyQ),
$$
where $ B_\varepsilon(\inftyQ)$ is a small open ball of hyperbolic 
volume $\varepsilon$ around the cusp $\inftyQ$. 
Since $0 \le \nu_\can \le \kappa_\circ \vol(X(\Gamma_{\!\!\circ})) \displaystyle{ \frac{ 
d}{g}} \,\mu_\Gamma$, we have
\begin{align*}
\int_{X(\Gamma)_{\varepsilon}}
-\log\|l\| ^2 \, \nu_\can  
&\le 
\int_{X(\Gamma)_{\varepsilon}}
\left[-\log\|l\| ^2\right]^+  \nu_\can\\
&\le \kappa_{\circ} \vol(X(\Gamma_{\!\!\circ}))  \frac{d}{g}\,
\int_{X(\Gamma)_{\varepsilon}}
\left[-\log\|l\| ^2\right]^+  \mu_{\Gamma}\\
&\le \kappa_{\circ}  \vol(X(\Gamma_{\!\!\circ})) \frac{d}{g}\,
\int_{X(\Gamma_{\!\!\circ}) \setminus B_\varepsilon(\inftyQ)}
\left[-\log\|l\| ^2\right]^+  \mu_{\Gamma_{\!\!\circ}}\\
& \le \frac{d}{g} \kappa_1(\varepsilon),
\end{align*}
where $[-\log\|l\| ^2]^+ =
(|-\log\|l\| ^2| -\log\|l\| ^2)/2$ is the positive
part of $-\log\|l\| ^2$ and $\kappa_1(\varepsilon)$ is a constant 
independent of $\Gamma$.   

It remains to bound the integral over the closure of $f^*
B_\varepsilon(\inftyQ)$. 
If $b$  denotes the width of the cusp $\inftyQ$, then we have 
\begin{align*}
 B_\varepsilon(\inftyQ) &= \{ z = x+i y \in \HH \,|\, 0<x \le b , \, y > b/\varepsilon \}
\end{align*}
and  for a cusp  $S_j$ ramified of order $b_j$ above $\inftyQ$
\begin{align*}
 B_{b_j\varepsilon}(S_j) &= \{ z = x+i y \in \HH \,|\, 0<x \le b b_j , \, y > b/\varepsilon \},
\end{align*}
therefore
\begin{align*}
f^* B_\varepsilon(\inftyQ) = \bigcup_{S_j} B_{b_j \varepsilon}(S_j).
\end{align*}
We may assume that the order of vanishing of $l$ at $\inftyQ$ is $1$.

For a local coordinate $t$ for the cusp $\inftyQ$ we shall write 
 $-\log\|l\|^2(t) = -\log|t|^2 + h(t)$, where $|h(t)|$ has
at most log-log growth.  As before we find  a constant 
 $\kappa_2(\varepsilon)$, which is independent of $\Gamma$, such that 
\begin{align*}
  \int_{\overline{B_{b_j\varepsilon}(S_j)}} 
    h(t) \nu_\can 
&\le \kappa_{\circ} \vol(X(\Gamma_{\!\!\circ})) \frac{d}{g}
  \int_{\overline{B_{b_j\varepsilon}(S_j)}}
  |h(t)| \mu_{\Gamma}\\
  &\le \kappa_{\circ} \vol(X(\Gamma_{\!\!\circ}))\frac{d}{g} \frac{b_j}{d}
  \,\int_{\overline{B_\varepsilon(\inftyQ)}}
  |h(t)| \mu_{\Gamma_{\!\!\circ}}\\
  &\le \frac{b_j}{g} \kappa_2(\varepsilon). 
\end{align*}
 
For the remaining term we recall that  the 
local coordinate at the  cusp $\inftyQ$ on the modular curve $X(\Gamma_{\!\!\circ})$ 
is given by $\exp(2 \pi i z/b)$ with $z=x +iy\in \HH$ and
$b\in \RR$. 
Therefore $-\log|t|^2 = 4 \pi y/b$ and 
\begin{align*}
  \int\limits_{\overline{B_{b_j\varepsilon}(S_j)}} -\log|t|^2
  \,\nu_\can 
  &=\int\limits_{\overline{B_{b_j\varepsilon}(S_j) \setminus
      B_{\varepsilon}(S_j)}} -\log|t|^2
 \, \nu_\can
  +\int\limits_{\overline{B_{\varepsilon}(S_j)}}  -\log|t|^2 \, \nu_\can\\
&\le   \kappa_\circ \frac{1}{g}
\iint\limits_{\substack{0<x<b_j b  \\\frac{b}{\varepsilon} < y < \frac{b b_j}{\varepsilon}}} 
\frac{4   \pi}{b} y  \, \frac{dxdy}{ y^2} 
  +\int\limits_{\overline{B_{\varepsilon}(S_j)}} -\log|t|^2\, \nu_\can\\
&=4 \pi  \kappa_\circ  \frac{b_j}{g} \log(b_j)
 +\int\limits_{\overline{B_{\varepsilon}(S_j)}} -\log|t|^2 \, \nu_\can.
\end{align*}

But now with the same considerations as in \cite{JK2}, p.1274, we get 
on 
\begin{align*}
\overline{B_\varepsilon(S_j)} =\{ z = x+i y \in \HH \,|\,
 0<x \le b b_j , \, y >b b_j/\varepsilon \}
\end{align*}
the bound
\begin{align}\label{eq:jokr}
g \,\nu_\can  &= y^2 \sum_{r=1}^g |f_r(z)|^2  \,\frac{dxdy}{y^2} \notag \\
&\le y^2 e^{- \frac{4\pi y}{b b_j}}  \max_{y = \frac{b b_j}{\varepsilon}} 
\sum_{r=1}^g |f_r(z)/ e^{- \frac{4\pi y}{b b_j}} |^2  \,\frac{dxdy}{y^2} \notag\\
&\le y^2 e^{-\frac{4\pi y}{b b_j}}  e^{\frac{4\pi}{\varepsilon_{\phantom{j}}} } \left(\frac{ 
 \varepsilon}{b b_j}\right)^2
\,\kappa_\circ \,\frac{dxdy}{y^2}.  \nne
\end{align}



Then \eqref{eq:jokr} yields 
\begin{align*}
  \int_{\overline{B_\varepsilon(S_j)}} -\log|t|^2 \nu_\can
&\le \frac{4 \pi}{b} \frac{\kappa_\circ}{g}  
\iint\limits_{\substack{ 0<x\le bb_j\\  \frac{b_j}{\varepsilon} <y }}  
 y e^{-\frac{4\pi y}{b b_j}}  e^{\frac{4\pi}{\varepsilon_{\phantom{j}} }}
 \left(\frac{\varepsilon}{b b_j} \right)^2
 \,dxdy\\
&=\frac{4 \pi}{b} \frac{\kappa_\circ}{g}   
\frac{b_j  b \varepsilon ( \varepsilon + 4 \pi)}{16  \pi^2}\\  
&  \le \frac{b_j}{g}\kappa_3(\varepsilon),
\end{align*} 
where $\kappa_3(\varepsilon)$  is again a constant 
independent of $\Gamma$.

Now for  the general case we take into account that the sum $\sum_{S_j} b_j$ of the ramification 
indices above the cusps equals $d$  and thus we derive 
\begin{align*}
\int_{\overline{f^* B_\varepsilon(\inftyQ)}}
-\log\|l\|^2 \nu_\can \le 
  4 \pi \kappa_\circ \sum_{S_j}  \ord_{f(S_j)}(l) \frac{b_j \log(b_j)}{g} + 
 \frac{d}{g}  (\kappa_2(\varepsilon) +\deg(L) \kappa_3(\varepsilon)).
\end{align*}
Setting $\kappa= \kappa_1(\varepsilon) +\kappa_2(\varepsilon) +
\deg(L) \kappa_3(\varepsilon)$ we deduce the claimed upper bound in
\eqref{bound:smoothint2}.

\hfill $\Box$

\section{Intersection numbers of hermitian line bundles}

In this section we  give an overview of Arakelov theory for arithmetic
surfaces, see e.g. \cite{Ga}, \cite{Kue},
\cite{So}.

\nn {\bf Notation.} \label{def:arithsurf} Let $K$ be a number field,
$\mathcal{O}_K$ its ring of integers and $\Sigma$ the set of complex
embeddings of $K$ in $\CC$.  An arithmetic surface
$\mathcal{X}$ is a regular scheme of dimension $2$ together with a
projective flat morphism $f: \mathcal{X} \longrightarrow \Spec
\mathcal{O}_K$, where $\mathcal{O}_K$ is the ring of integers of $K$.
Moreover, we assume that the generic fiber $X_K$ of $f$ is
geometrically irreducible and has genus greater than one.  

If $X_K$ is a smooth, projective curve
defined over $K$, then a regular model for
$X_K$ over $\Spec \mathcal{O}_K$ defines an arithmetic surface $\mathcal{X}$.  

We let $\calX_\infty$ be the set of
complex points   of the   $\Spec\,
\ZZ$-scheme  induced by $\mathcal X$, i.e.  $ \mathcal{X}_\infty =
\prod_{\sigma \in \Sigma} X_\sigma(\CC)$.   Note that the complex
conjugation $F_\infty$ acts on $\mathcal{X}_\infty$.  Finally, by
abuse of notation we set $$
\int_{\mathcal{X}_\infty} := \sum_{\sigma
  \in \Sigma} \int_{X_\sigma(\CC)}. $$

\nn {\bf Arithmetic Picard group.} We call a line bundle ${\mathcal L}$ on ${\mathcal X}$ equipped
with a hermitian metric $h$ a {\it hermitian line bundle} and denote
it by $\overline{{\mathcal L}}=({\mathcal L},h)$. As usual, 
$\c1(\overline \calL)$ is its first Chern form.  
 Let $\overline{{\mathcal
     L}},\overline{{\mathcal M}}$ be two hermitian line bundles on
 ${\mathcal X}$ and $l,m$ (resp.)  be non-trivial, global sections,
 whose induced divisors on ${\mathcal X}_\infty$ have no points in
 common. Then, the \emph{arithmetic intersection number
   $\overline{{\mathcal L}}.\overline{{\mathcal M}}$ of
   $\overline{{\mathcal L}}$ and $\overline{{\mathcal M}}$} is given
 by
  \begin{align} \label{eq:vaz}
    \overline{{\mathcal L}}.\overline{{\mathcal M}}:= (l.m)_{\rm
      fin}+( l.m )_\infty; \nne
  \end{align}
  here the contribution at the finite primes $(l.m)_{\rm fin}$ is
  defined by Serre's $\text{Tor}$-formula and the contribution at the
  infinite primes is given by a certain star product $( l.m )_\infty$
  (cf. the special cases \eqref{eq:star1}, \eqref{eq:star2} and
  \eqref{eq:star3} below).  Two hermitian line bundles
  $\overline{{\mathcal L}}, \overline{{\mathcal M}}$ on ${\mathcal X}$
  are {\it isomorphic}, if $ \overline{{\mathcal L}} \otimes
  \overline{{\mathcal M}}^{-1} \cong ({\mathcal O}_{\mathcal X} ,
  |\cdot |).  $ The \emph{arithmetic Picard group}, denoted by
  $\widehat{\textrm {Pic}}({\mathcal X})$, is the group of isomorphy
  classes of hermitian line bundles $\overline{{\mathcal L}}$ on
  ${\mathcal X}$, the group structure being given by the tensor
  product. The main theorem in Arakelov theory is that the
  formula \eqref{eq:vaz} induces a bilinear, symmetric pairing
\begin{eqnarray*}
\Pica({\mathcal X}) \times
 \Pica({\mathcal X})
 &\longrightarrow& \RR.
\end{eqnarray*}

Finally, $\Pica^{0}(\mathcal{X})_\QQ \subset \Pica(\mathcal{X})
\otimes_\ZZ \QQ$ denotes the subgroup generated by those hermitian
line bundles $\overline{\mathcal{L}}= (\mathcal{L}, \|\cdot \|)$
satisfying $\deg (\mathcal{L}|_{ \mathcal{C}^{(\mathfrak{p})}_{l}})=0$
for all irreducible components $\mathcal{C}^{(\mathfrak{p})}_{l}$ of
the fibers $f^{-1}(\mathfrak{p})$ above $\mathfrak{p} \in \Spec
\mathcal{O}_K$, and $\text{c}_1(\overline{\mathcal{L}}) =0$.

\nn {\bf Arithmetic Chow groups.}\label{nn:acg} Instead of the
arithmetic Picard groups $\widehat{{\rm Pic}}({\mathcal X})$ we can
also consider the \emph{arithmetic Chow groups} $\cha({\mathcal X})$.
The elements of these arithmetic Chow groups are equivalence classes of
arithmetic cycles represented by pairs $(D,g_D)$, where $D$ is a
divisor on $\calX$ and $g_D$ is a Green function for $D_\infty=
\prod_{\sigma \in \Sigma} D_\sigma(\CC)$.  If $\overline{\calL}$ is a
hermitian line bundle on $\mathcal X$, then the first arithmetic Chern class
$\ca1(\overline{\calL})\in \cha^1({\mathcal X})$ can be represented by
any of the pairs $(\dv(s), - \log\|s\|^2)$ where $s$ is a section of
$\calL$.  The assignment $\overline{\calL} \mapsto
\ca1(\overline{\calL})$ induces a morphism $\widehat{{\rm
    Pic}}({\mathcal X}) \to \cha^1({\mathcal X}) $ that is compatible
with the arithmetic intersection pairings in both groups. In
particular we have the formulas
\begin{align}
&(\dv(s),-\log\|s\|^2).(0,g) = \frac{1}{2} \int_{\calX_\infty} g \operatorname{c}_1
(\overline\calL),\label{eq:star1} \nne \\
&(0,g_1).(0,g_2) = \frac{1}{2} \int_{\calX_\infty} g_1 \; \operatorname{dd}^c g_2 
= \frac{1}{2} \int_{\calX_\infty} g_2 \; \operatorname{dd}^c g_1 
 \label{eq:star3} \nne,\\ 
&(\dv(s),-\log\|s\|^2).(D,0) = (\dv(s),D)_\fin, \label{eq:star2} \nne
\end{align}
which are just reformulations of
particular cases of \eqref{eq:vaz}.

\nn {\bf Admissible metrics.} \label{def:admissible} Let $\calX$ be an arithmetic surface.
Then $\calX_\infty = \prod_{\sigma \in \Sigma} X_\sigma(\CC)$ is a
finite union of compact Riemann surfaces. By abuse of notation we call
a $(1,1)$-form $\nu$ on $\calX_\infty$ such that $\nu= \prod_{\sigma
  \in \Sigma} \nu_\sigma$, where each $\nu_\sigma$ is a volume form
on $X_\sigma(\CC)$, also a volume form on $\calX_\infty$.  A
hermitian line bundle $\overline \calL$ is called
\emph{$\nu$-admissible} if $\operatorname{c}_1(\overline \calL)=
\deg(\calL) \nu$.

For a divisor $D = \sum n_P P$ on $\calX$ we define $
\overline\calO(D)_\nu=  (\calO(D),\|\cdot\|_\nu)
$, where the metric on the line bundle $\calO(D_\sigma)$ on
$X_\sigma(\CC)$ is such that $-\log \|1_{D_\sigma}\|^2_\nu= \sum n_p
g_{\nu_\sigma}(z,P)$, where $g_{\nu_\sigma}(z,P)$ is the unique
normalized Green function for $P_\sigma(\CC)$ associated with
$\nu_\sigma$ (see definition \ref{def:greenfunc}).  Observe that
$\overline\calO(D)_\nu$ is a $\nu$-admissible line bundle.
If $D$ is a fibral divisor we also write $\calO(D)$ for the hermitian line bundle $\calO(D) \otimes (\calO_\calX, |\cdot |)$.

%
%

\nn {\bf Logarithmically singular, hermitian line bundles.}  \label{def:hypadmissible} Let
$\calX$ be an arithmetic surface such that $\calX_\infty \cong
\prod_{\sigma \in \Sigma} X(\Gamma_\sigma)$ with fuchsian groups $\Gamma_\sigma$.  By abuse of notation we
write $\mu_\Gamma = \prod_{\sigma \in \Sigma} \mu_{\Gamma_\sigma}$.
It is natural to consider $\mu_\Gamma$-admissible line bundles
$\overline\calL$ on $\calX$, i.e., those metrized line bundles
equipped with a hermitian metric satisfying
$\operatorname{c}_1(\overline\calL)= \deg(\calL) \mu_\Gamma$.  Such
line bundles are examples of logarithmically singular, hermitian line
bundles in the sense of \cite{Kue}. It is shown in \cite{Kue} how to
modify the arithmetic intersection product in order to obtain an
arithmetic intersection product for such singular line bundles.  Here
we only need that for the particular cases discussed in subsection \ref{nn:acg}
the same formulae \eqref{eq:star1}, \eqref{eq:star3} and
\eqref{eq:star2} hold. We also write $\Pica(\calX)$ for the arithmetic Picard group of logarithmically singular hermitian line bundles.

\nn {\bf Theorem.} \label{thm:compare-metrics} \emph{Let $\nu$ be a
  volume form on $\calX_\infty$ and let $\mu$ either be another volume
  form or  a hyperbolic volume form.
  Then for any $\mu$-admissible hermitian
  line bundle $\overline\calL = (\calO(D),\|\cdot\|)$ of degree $d$ we
  have the equality
\begin{align*}
\overline{\calO}(D)_\nu ^2 
= \overline{\calL}^2 +  d\int_{\calX_\infty}
\log\|l\|^2 \nu + \frac{d^2}{2}c_{\mu,\nu},
\end{align*}
where $l$ is a section of $\calO(D)$ with divisor $D$ and 
$c_{\mu,\nu}$ is as in Lemma \ref{lem:gaac-hyp}.}

{\bf Proof.} Assume first that $\mu$ is a volume form.  In the
language of arithmetic Chow rings, which we will use for a moment, we
write
\begin{align*}
\widehat{\operatorname{c}}_1(\overline{\calO}(D)_\nu ) 
&= ( \dv(l), g_{\nu}(D,z)) \\
&= (\dv(l), -\log\|l\|^2) +  
(0, g_{\nu}(D,z) +\log\|l\|^2)) \\
&= \widehat{\operatorname{c}}_1(\overline{\calL}) +  
(0, g_{\nu}(D,z) - g_{\mu} (D,z)) +(0,g_{\mu}
(D,z)+\log\|l\|^2); 
\end{align*}
observe that the last term is in fact a vector of  constants. 
Then using formulae \eqref{eq:star1}, \eqref{eq:star3} and
\eqref{eq:star2}, 
we obtain   the equalities
\begin{align*}
\overline{\calO}(D)_\nu ^2 &= 
\overline{ \calL}^2 
+ \frac{d}{2} 
\int_{\calX_\infty}
 \big(g_{\nu}(D,z) - g_{\mu} (D,z) \big)\land \left(
   \mu(z) + \nu(z) \right)
+  d \left( g_{\mu}
(D,z)+\log\|l\|^2 \right)\\
&=\overline{ \calL}^2 
+   d \, a_{\mu,\nu}(D) + \frac{d^2}{2} c_{\mu,\nu} +  
d ( g_{\mu}
(D,z)+\log\|l\|^2)\\
&= \overline{\calL}^2 +  d \int_{\calX_\infty}
\log\|l\|^2  \nu + \frac{d^2}{2} c_{\mu,\nu}
\end{align*}
with the help of Lemma \ref{lem:intsecond}.
Now for $\mu_\Gamma$ we observe first that the same   proof
works. If $D$ has no support in the cusps, we just 
replace $g_\mu(D,z)$ by $g_{\mu_\Gamma}(D,z)$; else, if
the support of $D$ contains a cusp, we have to use $\tilde
g_{\mu_\Gamma}(D,z;w_0)$ and Lemma \ref{lem:intsecond-hyp} instead.
\hfill $\Box$

\nn {\bf Canonical metric.} 
If the genus of a Riemann surface $X$ is greater than one, then for each $\sigma \in \Sigma$ we
have on $X_\sigma(\CC)$ the \emph{canonical volume form}
\begin{align*}
\nu_\can^\sigma (z) = \frac{i}{2g} \sum_j |f_j^\sigma|^2 \dd z \land \dd \overline z,
\end{align*}
where $f_1^\sigma(z) \dd z$, ... $f_g^\sigma(z) \dd z$ is an
orthonormal basis of $H^0(X_\sigma(\CC),\Omega^1)$ equipped with the
natural scalar product. We write $\nu_\can$ for the induced volume
form on $\calX_\infty$ and for ease of notation we set 
$$
\overline\calO(D) = \overline\calO(D)_{\nu_\can}.
$$

\nn {\bf Adjunction formula.}  The following observation is due to Arakelov. There is a unique metric
  $\|\cdot \|_\Ar$ on $\omega_{\calX}$ such that for all sections $P$
  of $\calX$   the adjunction formula
\begin{align} \label{eq:adjunction}
\overline{\omega}_\Ar. \overline{\calO}(P)  + 
\overline{\calO}(P) ^2 = \log| \Delta_{K|\QQ}| \nne
\end{align}
holds. 
Moreover, $\overline{\omega}_\Ar =(\omega_{\calX},\|\cdot \|_\Ar)$ is a
$\nu_\can$-admissible line bundle.  Because of different conventions
we recall that
\begin{align} \label{def:omegaar}
\overline{\omega}_\Ar
=\overline{\omega}_{\calX/\calO_K,\Ar}\otimes f^*
\overline{\omega}_{\calO_K/\ZZ}, \nne
\end{align}
where $\omega_{\calO_K/\ZZ} =
\partial^{-1}_{K|\QQ}$ is equipped with the natural metric and
$\overline{\omega}_{\calX/\calO_K}$ with its residual metric (see e.g. \cite{MB}, eq. (1.2.1)). 
Observe that the adjunction  formula holds for arithmetic surfaces in the sense 
of subsection \ref{def:arithsurf} (see e.g. \cite{La}, p.101).

\nn {\bf Lemma.} \label{lem:wlee}\emph{Let $P$ be a section of 
$\calX \to \Spec \calO_K$ and  $\calF_P$ a fibral divisor
with the property 
\begin{align}\label{eq:def-calF}
\overline{\omega}_\Ar \otimes \overline{\calO}(P +\calF_P ) ^{2- 2g}
 \in \Pica^0(\mathcal{X})_\QQ.\nne
\end{align}
Then we have
\begin{align*}
\overline{\omega}_\Ar^2  =
- \operatorname{ht}_{NT}\left(\omega \otimes \calO(P)^{2- 2g}\right) +
 (4g -4)  \log| \Delta_{K|\QQ}|- 4g (g-1)\, 
\overline{\calO}(P) ^2 + (2g-2)^2 \calO(\calF_P)^2, 
\end{align*}  
where $\operatorname{ht}_{NT}$ denotes the Neron-Tate height on the Picard
group $\operatorname{Pic}^0(X)$
of $X$. In particular,
\begin{align*}
\overline{\omega}_\Ar^2  \le (4g -4)  \log| \Delta_{K|\QQ}|- 4g (g-1)\, 
\overline{\calO}(P) ^2. 
\end{align*}  
}

{\bf Proof.}  These formulae are well-known for semi-stable arithmetic
surfaces (see e.g. \cite{Fa}, p.410); however the essential
ingredients of their proof, namely the Faltings-Hriljac theorem (see
\cite{Hr} Theorem 3.1 and \cite{Hr} Proposition 3.3) and the
adjunction formula \eqref{eq:adjunction}, 
hold more generally for arithmetic surfaces in the sense of
subsection \ref{def:arithsurf}. 

 For convenience we recall the proof:
 There exists a fibral divisor $\calF_p$ which satisfies \eqref{eq:def-calF}.
Then, because of the Faltings-Hriljac theorem
, the definition of
 $\Pica^0(\mathcal{X})_\QQ$ and the adjunction formula
 \eqref{eq:adjunction}, we have
 \begin{align*}
- \operatorname{ht}_{NT}&\left(\omega\otimes \calO(P) ^{2- 2g} \right)\\
&= \left(\overline{\omega}_\Ar \otimes \overline{\calO}(P +\calF_P)^{2-
   2g}\right).\left( 
 \overline{\omega}_\Ar \otimes \overline{\calO}(P+\calF_P) ^{2- 2g}
   \right)\\ 
 &= \overline{\omega}_\Ar^2 + (4-4g) 
 \overline{\omega}_\Ar . \overline{\calO}(P)  +  (2g-2)^2 
 \overline{\calO}(P) ^2 - (2-2g)^2 \calO(\calF_P)^2 \\
 &=\overline{\omega}_\Ar^2  + \left( (4g-4) +  (2g-2)^2
 \right)\overline{\calO}(P) ^2  
 +  (4 -4g) \log| \Delta_{K|\QQ}|- (2-2g)^2 \calO(\calF_P)^2.
 \end{align*}
 Hence 
 \begin{align*}
   \overline{\omega}_\Ar^2 =- \operatorname{ht}_{NT}\left(\omega\otimes
   \calO(P) ^{2- 2g} \right) +(4g -4) \log|
     \Delta_{K|\QQ}|- 4g (g-1)\, \overline{\calO}(P) ^2
   + (2g-2)^2 \calO(\calF_P)^2, 
 \end{align*}
 and using the facts $\operatorname{ht}_{NT}\left(\omega\otimes
   \calO(P) ^{2- 2g} \right)\ge0$ and $\calO(\calF_P)^2 \le 0$ we
 derive the desired inequality. \hfill $\Box$

Note that $\calF_P$ as above 
has support in the fibers of $\calX$ above bad primes in the sense 
of subsection \ref{assumptions}.

\section{A new formula for $ \overline{\omega}_\Ar ^2$}
Here we consider  morphisms 
of arithmetic surfaces $\belyi : \calX \to \calY$ and relate $\overline{\omega}_\Ar ^2$ 
to arithmetic intersection numbers on $\calY$.
Clearly such a formula must  contain  terms
coming from the difference of the metrics on the line bundles in question
as well as some contributions coming from the primes of bad reduction 
of $\calX$.

\nn {\bf Proposition.}\label{prop:trace} \emph{ Let $\belyi: \calX \to
  \calY$ be a morphism of arithmetic surfaces.  Write $d$ for the
  degree of the induced morphism $\belyi: \calX_K \to \calY_K$. Let
  $P$ be $K$-rat\-ional point on $\calY_K$ and write $\belyi^*P =
  \sum_j b_j S_j$. We assume that all $S_j$ are $K$-rational points.
  We also denote by $P$, resp. $S_j$, the horizontal divisor induced
  by the divisor $P$, resp. $S_j$.  Let $\nu$ be a volume form on
  $\calX_\infty$ and let $\mu$ either be another volume form or a
  hyperbolic volume form.
 Finally, for $\overline{\calL} = (\calO(P),\|\cdot\|) \in
  \Pica(\calY)$, we assume that 
  $\belyi^* \overline{\calL}$ is a $\mu$-admissible hermitian line
  bundle. Under the above assumptions we have
  \begin{align} \label{eq:trace} \sum_j b_j\,
    \overline{\calO}(S_j)_\nu ^2\, = &\,\overline\calL^2
    +\int_{\calX_\infty} \log\|l\|^2 \nu + \frac{d}{2} c_{\mu,\nu}
    + \sum_j b_j\,  \calO(\calG_j)^2  \notag\\
    &- \sum_j b_j \operatorname{ht}_{NT}\big(
    \calO(S_j) \otimes \belyi^*\calL^{\otimes
      -1/d}\big),\nne
\end{align}
where $l$ is the  section of $\belyi^*\calO(P)$ whose divisor equals the
divisor $\belyi^*P$
 and for each cusp  $S_j$ the
vertical divisor  $\calG_j$ is defined in \eqref{eq:def-calg} below.}

{\bf Proof.} For each $S_j$ we consider a $\mu$-admissible, hermitian
line bundle $\calL_j=( \calO(S_j), \|\cdot\|)$ together with a fibral
divisor $\calG_j$ such that
\begin{align}\label{eq:def-calg}
  \left(\overline{\calL}_j \otimes \calO(\calG_j)\right) \otimes
  \left(\belyi^*\overline{\calL}\right)^{\otimes -1/d} \in
  \Pica^0(\calX)_\QQ.\nne
\end{align} 
Let $l_j$ be a section of $\calL_j$ with divisor $S_j$.  Then,
because of Theorem \ref{thm:compare-metrics} we have the equality
\begin{align}\label{eq:fslj}
  \overline{\calO}(S_j)_\nu ^2 = \overline{\calL}_j^2 +
  \int_{\calX_\infty} \log\|l_j\|^2 \nu + \frac{1}{2}c_{\mu,\nu}\nne .
\end{align}
Furthermore, we have
\begin{align*}
- \sum_j& b_j \operatorname{ht}_{NT}\big(
    \calO(S_j) \otimes \belyi^*\calL^{\otimes
      -1/d}\big) \notag\\
&= \sum_j b_j\, \left(\left(\overline{\calL}_j \otimes
  \calO(\calG_j)\right) \otimes
\left(\belyi^*\overline{\calL}\right)^{\otimes -1/d}\right)^2 \notag\\
&=\sum_j b_j\, \left( \overline\calL_j^2 + 2 \calL_j.\calO(\calG_j) + \calO(\calG_j)^2 
 -\frac{2}{d}  \,  
\left(\overline{\calL}_j
  \otimes \calO(\calG_j)\right).
\belyi^*\overline{\calL} + \frac{1}{d}\overline\calL^2 \right) \notag \\
&=
\sum_j b_j\, \left( \overline\calL_j^2 + 2 \calL_j.\calO(\calG_j) + \calO(\calG_j)^2 
 -\frac{2}{d} \,\calO(\calG_j).
\belyi^*\overline{\calL} \right) 
-\overline\calL^2 
+ \sum_j  b_j  \log\|l_j\|^2  - \log\|l\|^2  \notag \\
&=
\sum_j b_j\,  \overline\calL_j^2  
+ \sum_j  b_j  \log\|l_j\|^2
-\overline\calL^2 
  - \log\|l\|^2  
- \sum_j b_j\,  \calO(\calG_j)^2 ;\notag 
\end{align*}
observe that in the last equalities we used the relations
\begin{align*}
\sum_j b_j \,\overline{\calL}_j.
\belyi^*\overline{\calL}&=
\sum_j b_j ( S_j , -\log\|l_j\|^2 ). \belyi^*( P,  -\log\|l\|^2) \\
&=
\left( ( \sum_j b_j\,  S_j , - \log\|l\|^2 ) +
(0, - \sum_j b_j \, 
\log\|l_j\|^2 +\log\|l\|^2)\right). \belyi^* (P, -\log\|l\|^2) \\
&= \belyi_* ( \sum_j b_j\,  S_j , - \log\|l\|^2 ). (P, -\log\|l\|^2)
- d\,
 \left(\sum_j b_j  \log\|l_j\| -\log\|l\|  \right)\\
&= d \,\overline\calL^2 - d\,
 \left(\sum_j b_j  \log\|l_j\| -\log\|l\|  \right)
\end{align*}
and
\begin{align*}
\calO(\calG_j) . \frac{\belyi^*\calL}{d}& = \calL_j.\calO(\calG_j) + \calO(\calG_j)^2.
\end{align*}
Now  taking the weighted sum of \eqref{eq:fslj} 
the claim follows easily by recollecting the terms in 
the above equality for the Neron-Tate height in question.
\phantom{=}\hfill $\Box$

\nn {\bf Theorem.} \label{thm:keyformula} \emph{ 
Let $\belyi: \calX \to \calY$ be a morphism of arithmetic surfaces.
Write $d$ for the degree of the induced morphism  $\belyi: \calX_K \to
  \calY_K$. Let $P$ be $K$-rat\-ional
  point on $\calY_K$ and write $\belyi^*P = \sum_j b_j S_j$. We assume
that all $S_j$ are $K$-rational points.   
We also denote  by $P$, resp. $S_j$, the horizontal divisor induced
  by the divisor $P$, resp. $S_j$.
  Let $\mu$ 
either be a volume
  form or  a hyperbolic volume form.
Then for any $\overline{\calL} =
  (\calO(P),\|\cdot\|) \in \Pica(\calY)$ such that $\belyi^*
  \overline{\calL}$ is a $\mu$-admissible hermitian line
 bundle\footnote{unless $\belyi$ is unramified, the metric on $\calL$ must be hyperbolic.}, then
   the equality
 \begin{align}\label{eq:omegaformula1}
 \overline{\omega}_\Ar ^2=
& (4g -4)  \log|
  \Delta_{K|\QQ}| 
- \frac{4 g (g-1)}{d} \left(\overline\calL^2
    +\int_{\calX_\infty} \log\|l\|^2 \nu_\can + \frac{d}{2} c_{\mu,\nu_\can}
  \big)\right)\notag \\
&+ \frac{1}{d} \sum_j b_j
\Big(4 g (g-1)  \operatorname{ht}_{NT}\big(
    \calO(S_j) \otimes \belyi^*\calL^{\otimes
      -1/d}\big) -\operatorname{ht}_{NT}\left(\omega \otimes 
\calO(S_j)^{2- 2g}\right) \notag\\
&-4 g (g-1)\, \calO(\calG_j)^2
+(2g-2)^2 \,\calO(\calF_j)^2\Big)\nne 
 \end{align}
holds. Here $l$ is the unique section of $\belyi^*\calO(P)$ whose divisor 
equals the 
divisor induced by the divisor $\belyi^*P$ and 
$\calF_j$, resp. $\calG_j$, are as defined in \eqref{eq:def-calF}, resp. 
\eqref{eq:def-calg}. }

{\bf Proof.}  We write $\calF_j$ for the vertical divisor determined
by $S_j$ by means of formula \eqref{eq:def-calF}.  Again we also write
$S_j$ for the section induced by $S_j$.  Then for each $S_j$ we have
\begin{align*}
 \overline{\omega}_\Ar ^2& -(4g -4)  \log|
  \Delta_{K|\QQ}| \\
\quad&= -\operatorname{ht}_{NT}\left(\omega \otimes 
\calO(S_j)^{2- 2g}\right)+(2g-2)^2 \,\calO(\calF_j)^2 - (4 g (g-1)) \overline{\calO}(S_j)^2
\end{align*}
because of Lemma \ref{lem:wlee}. 
We now add these equalities
for all  $K$-rational points $S_j$ weighted with the factor
$b_j/d$, which is determined by $\belyi^*P= \sum b_j S_j$, and obtain
by means of \eqref{eq:trace} the  equality
\begin{align*}
 \overline{\omega}_\Ar ^2=
& (4g -4)  \log|
  \Delta_{K|\QQ}| 
- \frac{4 g (g-1)}{d} \left(\overline\calL^2
    +\int_{\calX_\infty} \log\|l\|^2 \nu_\can + \frac{d}{2} c_{\mu,\nu_\can}
  \big)\right)\\
&+ \frac{1}{d} \sum_j b_j
\Big(4 g (g-1)  \operatorname{ht}_{NT}\big(
    \calO(S_j) \otimes \belyi^*\calL^{\otimes
      -1/d}\big) -\operatorname{ht}_{NT}\left(\omega \otimes 
\calO(S_j)^{2- 2g}\right) \\
&-4 g (g-1)\, \calO(\calG_j)^2
+(2g-2)^2 \,\calO(\calF_j)^2\Big).
\end{align*}
\hfill $\Box$

\nn {\bf Remarks.} (i) If   the arithmetic
surface is semi-stable in Theorem \ref{thm:keyformula}, then the condition that all $S_i$ are
$K$-rational points can be skipped. Indeed, if $L$ is a finite Galois
extension of $K$ so that all $S_i$ are defined over $L$, then the
theorem holds for a minimal desingularisation $\pi: \calX' \to
\calX\times \Spec \calO_L$, i.e., the arithmetic surface obtained from $\calX$ via base change to $L$.  But then
due the fact that the canonical sheaf satisfies $\omega_{\calX'}
=\pi^*\omega_{\calX} \otimes f^*\omega_{\calO_L/\calO_K}$ (see
e.g. \cite{La}, p.127) and that the right hand side in formula
\eqref{eq:omegaformula1} is invariant under the Galois group
$\operatorname{Gal}(L|K)$, the claim follows by dividing the
resulting formula by $[L:K]$.

(ii) In many cases the arithmetic surface $\calX$
can be chosen to be the minimal regular model, but we stress
the fact that in general the arithmetic surface $\calX$ does not have to be  the
minimal regular model.  If $\calX$ is not the minimal regular model
$\calX_{\textrm{min}}$ of $X_K$, then there exists a morphism $\pi:
\calX \to \calX_{\textrm{min}}$ and a vertical divisor $\calW$ on
$\calX$ whose support is in those fibers which contain a
$(-1)$-curve, such that $ \pi^*\omega_{\calX_{\textrm{min}}}
=\omega_{\calX} \otimes \calO(\calW)$.  By means of Theorem
\ref{thm:keyformula} we then have
$$
\overline{\omega}_{\calX_{\textrm{min}},\Ar}^2 =
  \pi^*\overline{\omega}_{
\calX_{\textrm{min}},\Ar}^2
=\overline{\omega}_{\calX,\Ar}^2 + 2 
\omega_{\calX}.\calO(\calW) +\calO(\calW)^2
= \overline{\omega}_{\calX,\Ar}^2 + 
\sum_{\frakp\, \textrm{bad}} b_\frakp\, \log
 \Nm(\frakp),
 $$
with some uniquely determined coefficients $b_\frakp \in \QQ$.

(iii) Recall that a prime $\frakp$ is said to be bad if the fiber of $\calX$
above $\frakp$ is reducible. Then obviously the contributions of the
fibral intersections to \eqref{eq:omegaformula1} have only support at
the bad primes. 

\nn {\bf Definition.} Keep the notation of Theorem \ref{thm:keyformula}.

(i) We call  
\begin{align}\label{eq:def-ap}
 \sum_{\frakp \,\,{\rm bad}} a_\frakp
  \log\Nm(\frakp)
&= - \frac{2 g 
}{d} \sum_j b_j\,\calO( \calG_j)^2 +
\frac{
2g-2}{d} \sum_j b_j\, \calO(\calF_j)^2  \nne 
\end{align}
the \emph{geometric contributions}.

(ii) We call
\begin{align}\label{eq:def-abelyi}
a_\belyi&= 
- \frac{2 g}{d} \left(\overline\calL^2
    +\int_{\calX_\infty} \log\|l\|^2 \nu_\can + \frac{d}{2} c_{\mu,\nu_\can}
  \big)\right)\nne
\end{align}
the \emph{analytic contributions}.

Observe  that because
of Lemma \ref{lem:lowbound}, resp. Lemma \ref{lem:bound-c-hyp} we have
$a_\belyi \le b_\belyi$, where
\begin{align}\label{eq:def-bbelyi}
b_\belyi&= 
- \frac{2 g}{d} \left(\overline\calL^2
    +\int_{\calX_\infty} \log\|l\|^2 \nu_\can \nne
  \big)\right).
\end{align}

\nn {\bf Corollary.} \label{cor:keyformula}
\emph{If in addition to the assumptions of Theorem \ref{thm:keyformula}
all the divisors $d S_j -\belyi^*P$ are torsion divisors, i.e., if
$\operatorname{ht}_{NT}\big( \calO(S_j) \otimes \belyi^*\calL^{\otimes
  -1/d}\big)=0$, then
 \begin{align}\label{eq:omega-inequality} 
 \overline{\omega}_\Ar ^2\le
& (2g -2) \left(   \log|
  \Delta_{K|\QQ}|^2 
+b_\belyi + \sum_{\frakp \,\,{\rm bad}} a_\frakp
  \log\Nm(\frakp) \right). \nne
 \end{align}}

{\bf Proof.} This
upper bound for $\overline{\omega}_{\calX,\Ar}^2$ holds because 
the Neron-Tate height defines a positive
definite quadratic form on $\operatorname{Pic}^0(X)_\QQ$.
\hfill$\Box$

\section{The analytic contributions}

The analytic contributions \eqref{eq:def-bbelyi}
in Theorem \ref{thm:keyformula} can be bounded whenever 
$X$ varies in a discrete family of morphisms to a fixed base curve $Y$ 
unramified outside a fixed set of points $\inftyQ, P_1,...,P_r$.

\nn {\bf Proposition.}\label{cor:bound-bbelyi}\emph{ Let $\calX \to \calY$ be a morphism of arithmetic surfaces
as in subsection \ref{assumptions}. Let $b_\belyi$ be as in
\eqref{eq:def-bbelyi}. Then there exist constants $\kappa_\circ,
\kappa$ independent of $\calX$ such that
\begin{align} \label{eq:bound-bbelyi} b_\belyi \le - \frac{2g}{d}
    \overline{\calL}^2 + [K:\QQ] \left(8 \pi \kappa_\circ \sum
      \frac{b_j \log b_j}{d} + 2 \kappa \right).\nne
\end{align}
In particular, there exist constants $\kappa_1, \kappa_2$ independent
of $\calX$ such that 
\begin{align*}
b_\belyi \le [K:\QQ] \left( \kappa_1 \log(\max_j \{b_j\}) +
  \kappa_2 \right).
\end{align*}  }
 
{\bf Proof.}   The condition 
on $Y \setminus \{\inftyQ, P_1,...,P_r\}$ being hyperbolic 
implies that for each embedding $\sigma:K \to \CC$ 
the Riemann surface $Y_\sigma(\CC)$ is
associated to a fuchsian group  $\Gamma_{\!\!\circ}^\sigma$ of the first kind.
Of course cusps and elliptic points are considered to be among $\inftyQ, P_1,...,P_r$. 
Moreover, $\belyi: \calX_\infty \to \calY_\infty$ is for each $\sigma$ a morphism 
of hyperbolic  curves   as in subsection \ref{CCModelle}.
Now the inequality \eqref{eq:bound-bbelyi} is a direct consequence
of Theorem \ref{thm:jk+}. The Hurwitz formula implies that $g/d$ is bounded,
thus the second claim is established.
\hfill $\Box$

\section{An upper bound for the geometric contribution}
Here we give general bounds for the quantity $ \sum_{\frakp\,\,\rm{
    bad}} a_\frakp \log \Nm(\frakp) $ defined by \eqref{eq:def-ap},
i.e., the geometric contributions in the formula of  Theorem \ref{thm:keyformula}.

For each bad prime $\frakp$ we let
\begin{align} \label{def:r_p}
\calX \times \overline
\FF_\frakp
= \sum_{j=1}^{r_\frakp} m_j C_j^{(\frakp)}\nne
\end{align}
be the decomposition into irreducible components
and set 
\begin{align}\label{def:u_p}
u_\frakp= \max_{i,j}
|C_i^{(\frakp)}.C_j^{(\frakp)}|, \qquad
l_\frakp= \min_{C_i^{(\frakp)}.C_j^{(\frakp)}\neq 0}
|C_i^{(\frakp)}.C_j^{(\frakp)}|.\nne
\end{align}
We further denote by $c_\frakp$ the 
connectivity of the dual graph of $\calX(\overline
\FF_\frakp)$, i.e. the 
minimal number
of intersection points needed to connect any two
irreducible components. We define
\begin{align}\label{eq:def-bp}
 b_\frakp = \left( 
 \sum_{k=1}^{c_\frakp} 
\left( \sum_{l=1}^{ k} 
 \left(\frac{u_\frakp}{l_\frakp}\right)^{l-1}  \right)^2
+ (r_\frakp - c_\frakp-1) \left( \sum_{l=1}^{ c_\frakp} 
 \left(\frac{u_\frakp}{l_\frakp}\right)^{l-1}  \right)^2 \right)  
 \frac{u_\frakp}{l_\frakp^2}  .\nne 
\end{align}

\nn {\bf Proposition.} \label{prop:gbound} \emph{
Let $\calG_j$ be as in \eqref{eq:def-calg}. Then  
its  component $\calG_j^{(\frakp)}$ in the fibre above a bad prime $\frakp$ satisfies
\begin{align*}
- \left(\calG_j^{(\frakp)}\right)^2 \le b_\frakp.
\end{align*}
}

{\bf Proof.}
After possibly
renumbering the irreducible components and adding rational multiples
of entire fibers, we may assume $0 \neq\calG_j^{(\frakp)}=
\sum_{k=2}^{r_\frakp} n_k C_k^{(\frakp)}$ with all $n_k\ge 0$ and
$n_1=0$.

 Let $W=\{\,C_j^{(\frakp)}\}$ be the set of  irreducible components 
of the fibre above $\frakp$ and set
\begin{align*}
U_0 &= \{ C_j^{(\frakp)} \in W \,|\, n_j =0 \,\} \\
V_0 &= W \setminus U_0.
\end{align*}
Then we define recursively 
\begin{align*}
U_{k+1} &= \{ C_j^{(\frakp)} \in V_k \,|\, \exists\, C_i^{(\frakp)} \in U_k \,\mbox{with}\, C_j^{(\frakp)}.C_i^{(\frakp)} > 0 \,\} \\
V_{k+1} &= V_k \setminus U_{k+1}.
\end{align*}
Since the fibre above $\frakp$ is connected, the subsets $U_k \subset W$ determine
a disjoint decomposition of $W$. In fact this decomposition has 
at most  $c_\frakp+1$ disjoint sets.

For all $ C_j^{(\frakp)} \in U_1$ we have 
$$n_j \le \frac{1}{l_\frakp}.$$
Indeed 
there exists a $C_l^{(\frakp)} \in U_0$ with $
C_j^{(\frakp)}. C_l^{(\frakp)}>0$ and using this component
 we obtain
 the  upper bound
\begin{align} \label{eq:bound-n2}
\deg \left(\calL|_{C_l^{(\frakp)}}\right) =
d \left(  S_j + \calG_j^{(\frakp)}\right) . C_l^{(\frakp)} = d  \left(  S_j. C_l^{(\frakp)}  +
\sum_{k=2}^{r_\frakp} n_k
C_k^{(\frakp)}. C_l^{(\frakp)}\right) \le d.\nne
\end{align}
Recall that $S_j. C_l^{(\frakp)} \ge 0$ and that by construction all the
summands $n_k C_k^{(\frakp)}. C_l^{(\frakp)}$ in \eqref{eq:bound-n2}
are non-negative. Thus we can derive the following bound for $n_j$:

For all $ C_j^{(\frakp)} \in U_2$ we have 
$$n_j \le  \frac{1}{l_\frakp}\left( 1 + \frac{u_\frakp}{l_\frakp}\right) .$$
Indeed there exists a $C_l^{(\frakp)} \in U_1$ with $
C_j^{(\frakp)}. C_l^{(\frakp)}>0$ and using this component we obtain the
upper bound
\begin{align*}
d n_j l_\frakp &\le d n_j  C_j^{(\frakp)}. C_l^{(\frakp)} \le 
 d S_j. C_l^{(\frakp)} + d \sum_{
\substack{k=2\\k \neq l}}
^{r_\frakp} n_k
C_k^{(\frakp)}. C_l^{(\frakp)} 
= \deg \left(\calL|_{C_l^{(\frakp)}}\right)
- d  n_l\, C_l^{(\frakp)}. C_l^{(\frakp)}\\
&\le d \left(1 + \frac{u_\frakp}{l_\frakp}\right),
\end{align*}
since $0 \le  \deg \left(\calL|_{C_l^{(\frakp)}}\right) \le d$ and
$  C_l^{(\frakp)}. C_l^{(\frakp)} <0$.

Repeating this procedure we get for the remaining coefficients
$n_j$ of $ C_j^{(\frakp)} \in U_k$ the bound
$$n_j \le  \frac{1}{l_\frakp} \sum_{l=1}^{k}  \left( \frac{u_\frakp}{l_\frakp}\right)^{l-1}.$$

So we obtain
\begin{align}\label{eq:betterbound}
  - \left(\calG_j^{(\frakp)}\right)^2 
&= - \sum_{j, k=2}^{r_\frakp} n_j
    n_k C_j^{(\frakp)}. C_k^{(\frakp)}  \notag\\
&\le -\sum_{j=2}^{r_\frakp} n_j^2
     C_j^{(\frakp)}. C_j^{(\frakp)}  \notag\\
&\le 
 \sum_{\substack{ U_k \subset W \\ U_k \neq U_0}} \# U_k \cdot  
\left( \sum_{l=1}^{k}  \left(\frac{u_\frakp}{l_\frakp}\right)^{l-1}  \right)^2 \frac{u_\frakp}{l_\frakp^2}  \notag \\
&\le b_\frakp.
\nne
\end{align}
For the last inequality we used $1 \le \# U_k$, $ \# U_{c_\frakp} \le r_\frakp - c_\frakp $
and $\frac{u_\frakp}{l_\frakp} \ge 1$.
 \hfill $\Box$

\nn {\bf Proposition.} \label{prop:g=f} \emph{Let $\calG_j$, resp. $\calF_j$, be 
as in \eqref{eq:def-calg}, resp. \eqref{eq:def-calF}.  If $\belyi: \calX \to \PP^1 $ is a
Galois cover, i.e. the extension of the function fields 
$K(\calX) \to  K(\PP^1_{\mathcal{O}_K}) $ is Galois with group $G$ and $\PP^1_{\mathcal{O}_K}$ is isomorphic to $\calX /G$,
then we have
\begin{align*}
  \calG_j^2 = \calF_j^2.
\end{align*}}

{\bf Proof.}
If $f: \calX \to \calY$ is a Galois
  cover of arithmetic surfaces, then there exists a line
bundle
$\calL$ on $\calY$ with:
\begin{align}\label{claim:galois}
\omega_{\calX}^{\deg(f)} \cong f^* \calL.\nne
\end{align}
Indeed, by assumption $\calY = \calX / G$ for a finite
group $G$. Let $s \in H^0(\calX, \omega_\calX)$, then also $s^\sigma
\in H^0(\calX, \omega_\calX)$ for all $\sigma \in G$ and therefore
$$
h= \prod_{\sigma \in G} s^\sigma  \in  H^0(\calX,
\omega_\calX^{\otimes |G|}).
$$
However, since $h^\sigma = h$ for all $\sigma
 \in G$, we have $\calO(\dv(h)) = f^* \calL$ for some line bundle
$\calL$ on $\calY$. Hence
$$
\omega_\calX^{\otimes|G|} \cong \calO(\dv(h)) \cong f^* \calL,
$$
from which we deduce the claim \eqref{claim:galois}  using the fact
$|G| = \deg(f)$.

Since
any line bundle on $\PP^1$ is of the form $\calO(1)^{\otimes e}
\otimes  \calO(\calF)$, where the fibral divisor $\calF$ is a finite
sum of entire fibers, we  have $\calF_j  = \calG_j$. 
\hfill $\Box$

\nn {\bf Theorem.} \label{thm:geom-bound}
\emph{
Let $a_\frakp$ be as in \eqref{eq:def-ap} and 
$b_\frakp$ as in \eqref{eq:def-bp}, then
for all primes of bad reduction we have
\begin{align*}
a_\frakp \le 2g\, b_\frakp.
\end{align*}
If in addition   
$\belyi: \calX \to \calY$ is a Galois cover and if $\calY =\PP^1$,
 then   the stronger inequality
\begin{align*}
a_\frakp \le 2\, b_\frakp
\end{align*}
holds.
}

{\bf Proof.}  In order to verify the first claim 
it suffices
to bound the terms involving $\calO(\calG_j)^2$,
since $\calO(\calF_j)^2 \le 0$ for all $j$.  
Upper bounds have been given in Proposition \ref{prop:gbound}.
Now summing  over all $j$  yields the first
claim.

In order to verify  the second claim we first observe that by 
Proposition \ref{prop:g=f} the quantity $a_\frakp$ only involves
$\calG_j$ and then we use Proposition \ref{prop:gbound} again.
\hfill $\Box$

\section{Explicit calculations for the modular curves $X_0(N)$ }
We first   give 
explicit formulas for the coefficients $a_p$ 
in \eqref{eq:def-ap} for the modular curves
$X_0(N)$. This in turn proves Theorem II from the 
introduction.

\nn {\bf Preliminaries.}  The curves $X_0(N)$ 
have a model over $\QQ$ and its complex valued
points correspond to the compact Riemann surfaces $\Gamma_0(N)
\setminus \left( \HH \cup \PP^1(\QQ) \right)$, where $\Gamma_0(N)=
\left\{ \left(\begin{smallmatrix} a & b \\ c & d\end{smallmatrix}
  \right) \in \Gamma(1) \,\left| \,c \equiv 0 \mod{N} \right.
\right\}$. 
From now on $N$ will be a square free integer  with at least two different prime factors
such that $(6,N)=1$, because 
then the cusps are also rational
points \cite{Og}. The Manin-Drinfeld theorem assures that any divisor 
of degree zero with support in the cusps is a torsion divisor. 
The index in $\Gamma(1)$ is given by the formula
\begin{align*}
d= [\Gamma(1):\Gamma_0(N)] = \prod_{p|N} (p+1).
\end{align*}
Recall that in this situation the genus of $X_0(N)$  is given by
\begin{align*} 
g_N= 1 + \frac{d}{12} - \frac{1}{4} \prod_{p|N} \left( 1 + \left(
    \frac{-1}{p}\right )\right) - \frac{1}{3} \prod_{p|N} \left( 1 + \left(
    \frac{-3}{p}\right )\right) - \frac{1}{2} \sigma_0(N),
\end{align*}
where $(\frac{\cdot}{\cdot})$ is the Legendre symbol and $\sigma_0(N)$ denotes
the number of divisors of $N$. Therefore   we have 
\begin{align}\label{eq:genusx0n}
\frac{12(g-1)}{d} = 1 + O(N^{-\varepsilon})\nne
\end{align}
for a small $\varepsilon>0$, since $N$ is square free.
\nn {\bf Minimal model.} The minimal model $\calX_0(N)$ over $\Spec \ZZ$ has been
determined by Deligne and Rapoport to be as follows:  The curve
$X_0(N)$ is smooth over $\ZZ[1/N]$. If $p|N$, then the scheme
$\calX_0(N) \times \ZZ/p\ZZ$ is reduced and singular over $\ZZ/p\ZZ$.
We write $N=pM = p q_1 \cdot\ldots \cdot q_\nu$, set $Q=
\prod_{i=1}^\nu (q_i +1)$ and define
\begin{align*}
u & = \begin{cases} 1 & \textrm{if } p \equiv 7  \textrm{ or } 11
  \pmod{12}  \textrm{ and all } q_i \equiv 1 \pmod{4}, \;\; i=1\ldots \nu\\
0 & \textrm{otherwise} 
\end{cases}\\
v & = \begin{cases} 1 & \textrm{if } p \equiv 5  \textrm{ or } 11
  \pmod{12}  \textrm{ and all } q_i \equiv 1 \pmod{3}, \;\; i=1\ldots \nu\\
0 & \textrm{otherwise.} 
\end{cases}
\end{align*}

The fiber $\calX_0(N) \times \ZZ/p\ZZ$ is the union of two copies
$C_0$, $C_\infty$ of $\calX_0(M) \times \ZZ/p\ZZ$ crossing
transversely in certain supersingular points and some chains of
projective lines connecting the remaining supersingular points with
those of the other copy. More precisely, if $u=1$, then there are
$2^\nu$ projective lines $F_i$, if $v=1$, then there are $2^\nu$ pairs
$G_i$, $H_i$ of projective lines and if $u=v=0$, there are no such
projective lines. All the projective lines have self-intersection
number $-2$.  The intersection behavior between these irreducible
components is given by the figure 
below.


\begin{picture}(0,0)%
\includegraphics{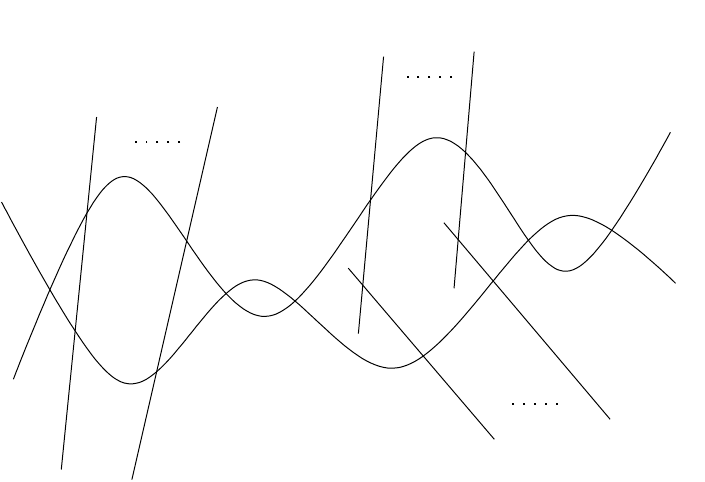}%
\end{picture}%
\setlength{\unitlength}{4144sp}%
\begingroup\makeatletter\ifx\SetFigFont\undefined%
\gdef\SetFigFont#1#2#3#4#5{%
  \reset@font\fontsize{#1}{#2pt}%
  \fontfamily{#3}\fontseries{#4}\fontshape{#5}%
  \selectfont}%
\fi\endgroup%
\begin{picture}(5355,3652)(4489,-4393)
\put(9790,-1734){\makebox(0,0)[lb]{\smash{{\SetFigFont{10}{12.0}{\rmdefault}{\mddefault}{\updefault}{\color[rgb]{0,0,0}$C_{\infty}$}%
}}}}
\put(9829,-2846){\makebox(0,0)[lb]{\smash{{\SetFigFont{10}{12.0}{\rmdefault}{\mddefault}{\updefault}{\color[rgb]{0,0,0}$C_0$}%
}}}}
\put(7373,-889){\makebox(0,0)[b]{\smash{{\SetFigFont{10}{12.0}{\rmdefault}{\mddefault}{\updefault}{\color[rgb]{0,0,0}$H_1$}%
}}}}
\put(8256,-4266){\makebox(0,0)[b]{\smash{{\SetFigFont{10}{12.0}{\rmdefault}{\mddefault}{\updefault}{\color[rgb]{0,0,0}$G_1$}%
}}}}
\put(9138,-4189){\makebox(0,0)[b]{\smash{{\SetFigFont{10}{12.0}{\rmdefault}{\mddefault}{\updefault}{\color[rgb]{0,0,0}$G_{2^{\nu}}$}%
}}}}
\put(5225,-1350){\makebox(0,0)[b]{\smash{{\SetFigFont{10}{12.0}{\rmdefault}{\mddefault}{\updefault}{\color[rgb]{0,0,0}$F_1$}%
}}}}
\put(6146,-1350){\makebox(0,0)[b]{\smash{{\SetFigFont{10}{12.0}{\rmdefault}{\mddefault}{\updefault}{\color[rgb]{0,0,0}$F_{2^{\nu}}$}%
}}}}
\put(8102,-889){\makebox(0,0)[b]{\smash{{\SetFigFont{10}{12.0}{\rmdefault}{\mddefault}{\updefault}{\color[rgb]{0,0,0}$H_{2^{\nu}}$}%
}}}}
\end{picture}%

In this figure  
all intersections are transversal, i.e., all
intersection multiplicities are equal to $1$. Finally we mention the 
formula
\begin{align*}
C_0.C_\infty = d \frac{p-1}{12(p+1)} - 2^\nu \left( \frac{u}{2} +
  \frac{v}{3}\right).
\end{align*}

The  natural morphism $\belyi: X_0(N) \to X(1)$ extends to a morphism 
 $\belyi: \calX_0(N) \to \calX(1)$. Its degree equals $d$ and 
the local degrees at the bad fibres are $d/(p+1)$, $pd /(p+1)$ or zero 
(resp.) at the irreducible components $C_\infty$, $C_0$ or at any of the other kind of fibres (resp.). 

There are two cusps $0$ and $\infty$ and  their associated 
horizontal divisors are disjoint.

\nn {\bf Lemma.} \label{lem:x0nG} \emph{
With  notation as in  \eqref{eq:def-calg}  
we have on 
$\calX_0(N)$ the formula
\begin{align*}
   \sum_{j} b_j \calO(\calG_j)^2
&=  -12 \sum_{p|N}  \frac{p}{p^2-1} \log p,
\end{align*}
}

{\bf Proof.}   
Recall that the global sections of the line bundle $\calL= \belyi^* \calO(1)$
are given by modular forms of weight $12$. Because of the
Manin-Drinfeld theorem there exists for any cusp $S_j$ an integral
modular form $g_j$ of some weight $12 e$ which only vanishes 
at $S_j$ on $X_0(N)$. 

We can consider the above equality
prime by prime.
 We put 
\begin{align*}
  \calD_0 &= C_0 + \frac{1}{2} \sum_{\nu}  F_\nu +
\frac{2}{3} \sum_\mu  G_\mu + \frac{1}{3} \sum_\rho  H_\rho;
\end{align*}
in the cases where some of the components $F_\nu, G_\mu, H_\rho$ do not exist
these are set to be the zero divisors.

 For a cusp $S_j$  we have
\begin{align*}
\dv(g_j^{1/e})\big|_{\calX_0(N)\otimes \ZZ/p \ZZ} 
 \equiv   d \,(S_j + \alpha_j \calD_0) \mod \calX_0(N)\otimes \ZZ/p \ZZ,
\end{align*}
where
\begin{align*}
\alpha_j = \frac{12}{d(p-1)} \left( 1-  (p+1)\,S_j.C_\infty\right) .
\end{align*} 
Indeed, the relations 
$
\deg \dv(g_j^{1/e})\big|_{F_\nu} =\deg \dv(g_j^{1/e})\big|_{G_\mu} =\deg\dv(g_j^{1/e})\big|_{H_\rho} =0
$
are easily checked and  also the remaining relation 
\begin{align*} 
\frac{d}{p+1}= \deg \dv(g_j^{1/e})\big|_{C_\infty} = 
d \,\left(S_j + \alpha_j \calD_0 \right).C_\infty 
=d \, \left(S_j.C_\infty + \alpha_j  d \frac{p-1}{12(p+1)}\right)
 \end{align*}
 holds. Set $ \calD_\infty = \calX_0(N)\times \ZZ/p\ZZ - \calD_0 =
 C_\infty + \frac{1}{2} \sum_{\nu} F_\nu + \frac{1}{3} \sum_\mu G_\mu
 + \frac{2}{3} \sum_\rho H_\rho$.  From $\calD_0.
 (\calD_0+\calD_\infty)=0$ we deduce
$$
\calD_0^2=-\calD_0.\calD_\infty = -\calD_0.C_\infty=- d
\frac{p-1}{12(p+1)}.
$$
Then as $ \calG_j^{(p)} = \alpha_j \,\calD_0$ and as $S_j.C_\infty \in \{0,1\}$, 
the last equality  implies 
\begin{align*}
( \calG_j^{(p)})^2   = 
 \alpha_j^2 \, \calD_0^2 = 
-  \frac{12}{d} \left( \frac{1}{p^2-1} +  \,S_j.C_\infty\right)  . 
\end{align*}
Now summing  over all cusps $S_j$ yields the identity
\begin{align*}
 \sum_{S_j} b_j \calO(\calG_j^{(p)})^2 
&=  - \frac{12}{d}
   \sum_{S_j} b_j \left( \frac{1}{p^2-1} + S_j.C_\infty   \right) \log p\\
&= - \frac{ 12}{ d} \left( \frac{d}{p^2-1}  +  \frac{d}{p+1}\right) \log p\\
&= -12 \left(\frac{p}{p^2-1}\right) \log p.
\end{align*}
Finally summing  over all the primes of bad reduction leads to the claim.
\hfill $\Box$ 

\nn {\bf Lemma.} \label{lem:x0nF} \emph{
With the notation as in \eqref{eq:def-calF} 
we have for the arithmetic surface 
$\calX_0(N)$ the formulae
\begin{align*}
 \sum_{p|N}  \sum_{S_j} b_j \calO(\calF_j^{(p)})^2
= -3 \sum_{p|N}  \frac{p+1}{p-1}\log p 
= -3 \log N  - \sum_{p|N}  \frac{6}{p-1}\log p. 
\end{align*}
}

{\bf Proof.} 
The proof of this formula is similar to the proof of Lemma
\ref{lem:x0nG}; we
refer to \cite{AU}, Proposition 4.2.1, p.63. \hfill $\Box$

\nn {\bf Proposition.} \label{prop:x0n} \emph{
With the notation as in \eqref{eq:def-ap}
we have for the arithmetic surface 
$\calX_0(N)$ the estimate
\begin{align*}
\sum_{p \,\,{\rm bad}} a_p \log p =- \frac{1}{2} \log N + O(N^{-\varepsilon}\log(N))
\end{align*}
}

{\bf Proof.} 
By means of the Lemmata \ref{lem:x0nG} and \ref{lem:x0nF} we have the equality
\begin{align*}
\sum_{p \,\,{\rm bad}} a_p \log p &=  - \frac{ 2 g}{d} \sum_{j} b_j \calO(\calG_j)^2
+ \frac{(2g-2)}{d} \sum_{S_j} b_j \calO(\calF_j)^2\\
&=  \frac{2}{d} \left(  
  \sum_{p|N} \left(12  g  \frac{p}{p^2-1}  -  (g-1) \frac{6}{p-1} \right)  \log p  - 3 (g-1) \log N \right)\\
&= - \frac{1}{2} \log N + O(N^{-\varepsilon}\log(N)).
\end{align*}
\hfill $\Box$

\nn {\bf Theorem.} 
 \emph{Let
$N$ be a square-free integer
   having at least two different prime factors and satisfying $(N,6)=1$.
Let $\calX_0(N)$ be the minimal regular model
   of the modular curve $X_0(N)$ and $g_N$ its genus. Then the
   arithmetic self-intersection number of its dualizing sheaf equipped
   with the Arakelov metric is bounded from above by
\begin{align}\label{eq:thmxon1}
\overline{\omega}_\Ar^2  
&\le    (16 \pi \kappa_\circ -1) g_N \log(N) +O(g_N), \nne
\end{align}
where $\kappa_\circ \in \RR$ is an absolute constant  independent of $N$.
}

{\bf Proof.} To ease notation we write $g$ for $g_N$.
The minimal model
$\belyi: \calX_0(N) \to \calX(1)$ together with its cusps $S_j$
satisfies the assumption of Theorem \ref{thm:keyformula}. We take 
$P$ to be the cusp $\infty \in X(1)(\QQ)$. We may assume that
$\overline\calL$ is $\mu_\Gamma$-admissible.

We now present bounds or  estimates for the positive terms on the right hand side of \eqref{eq:omegaformula1}. 

 Since $K=\QQ$, the  relative discriminant equals $1$, thus
the first term vanishes. 

By the genus formula
\eqref{eq:genusx0n} for $X_0(N)$ we deduce  
\begin{align*}
- \frac{4 g (g-1)}{d}  \overline\calL^2 = O(g)
\end{align*}
for the second  term. 
For the integral we use Theorem \ref{thm:jk+} and the fact that
the cusp $0$ has the maximal possible width $N$, i.e. $b_{\max}=N$, hence
\begin{align*}
\frac{4 g (g-1)}{d} \int_{X(\Gamma_0(N))} -\log\|l\| \, \nu_\can  &\le 
\frac{4 g (g-1)}{d} \left( 4 \pi \kappa_\circ \sum_{S_j}
 \frac{b_j \log(b_j)}{g} +  \frac{d}{g} \kappa\right)\\
&\le 
 (g-1)  16 \pi \kappa_\circ \log(N) + O(g).
\end{align*}
We can neglect the contribution involving the
 quantity $c_{\mu_\Gamma,\nu_\can}$, as   
this is a positive real number by Lemma \ref{lem:bound-c-hyp}. 
 
There is no positive contribution  of the terms involving Neron-Tate heights. 
This is due to the theorem of
 Manin-Drinfeld and the non-negativity of the Neron-Tate height. 
 
The remaining contributions  have been calculated in 
Proposition \ref{prop:x0n}
\begin{align*}
(2 g-2) \sum_{p \,\,{\rm bad}} a_p \log p =- g \log N + O(g)  .
\end{align*}
\hfill $\Box$

Since $\calX_0(N)$ is a semi-stable model of $X_0(N)$,
the arithmetic self-intersection number of its dualizing sheaf is positive by a theorem of Faltings, this 
fact immediately  implies  the following results.

\nn {\bf Corollary.} \emph{Keep the previous notation. Then
\begin{align*} 
\sum_j \frac{b_j}{d}\operatorname{ht}_{NT}\left(\omega \otimes 
\calO(S_j)^{2- 2g}\right) = O(g_N \log (N)) 
\end{align*}
and
\begin{align*} 
c_{\mu_\gamma, \nu_\can} = O( \log(N)/g_N).
\end{align*}
}
\hfill $\Box$

\section{Bounds for  Fermat curves $F_p$}
 
 We now consider the Fermat curves $F_n: x^n
+y^n =z^n$. The morphism $\belyi: F_N \to \PP^1$ given by $(x:y:z)
\mapsto (x^n:z^n)$ determines a Galois cover with Galois group $(\ZZ/n\ZZ)^2$.
Since $\belyi$ only has the three ramification points $0, 1, \infty$, it is
a Belyi morphism. The ramification indices are all equal to $n$.
In contrast to the previous examples, we proceed
with the identification $\PP^1 \setminus \{0,1,\infty\} = \Gamma(2)
\setminus \HH$.  The principal congruence subgroup $\Gamma(2)$ is a
free group on two generators $A,B$. Let 
$$\Gamma_n= \left\{ \gamma=
  A^{e_1}B^{f_1} \cdots A^{e_r}B^{f_r} \in \Gamma(2) \,\left| \, \sum
    e_i = \sum f_j = n \right.  \right\},$$
 then $\Gamma(2)/\Gamma_n
\cong (\ZZ/n\ZZ)^2$. Therefore $F_n = X(\Gamma_n)$ and $\belyi$ is
induced by the natural morphism $\Gamma_n \setminus \HH \to \Gamma(2)
\setminus \HH$. The cusps are defined over $\QQ(\zeta_n)$ and the
group of cuspidal divisors on $F_n$ with respect to the uniformization
given by $\Gamma_n$ modulo rational equivalence is a torsion subgroup;
in fact its structure is determined in \cite{Ro}.

Let $\calX$ be the desingularisation of the
  closure in $\PP^2_{\ZZ[\zeta_p]}$ of the Fermat curve $x^p+y^p=z^p$ 
with prime exponent $p$ as given in  \cite{Mc}. It is smooth outside the fibers over the set of primes
above $p$. If $\frakp$ is a prime above $p$, then it is shown in \cite{Mc} that
$\calX(\overline \FF_\frakp)$ has at most 
\begin{align} \label{eq:rp-fermat}
r_\frakp \le 4+ p(p-3)/2\nne
\end{align}
irreducible components $C^{(\frakp)}_l$, which are all isomorphic to $\PP^1$, and that
\begin{align} \label{eq:sp-fermat}
u_\frakp:=\max_{l,m}\left|C^{(\frakp)}_l.C^{(\frakp)}_m \right| \le p 
\mbox{ and } l_\frakp = \min_{C^{(\frakp)}_l.C^{(\frakp)}_m\neq 0}\left|C^{(\frakp)}_l.C^{(\frakp)}_m \right|=1.
\nne
\end{align} 
Furthermore, we have $c_\frakp=4$ for the connectivity of the dual graph of $\calX$, thus
\begin{align*}
b_\frakp &\le  \left( 
 \sum_{k=1}^{4} 
\left( \sum_{l=1}^{ k} 
 p^{l-1}  \right)^2
+  \frac{p(p-3)-2}{2}  \left( \sum_{l=1}^{4} 
 p^{l-1}  \right)^2 \right)  
 p \le \frac{ p^8}{2}. 
\end{align*}

\nn {\bf Theorem.} 
 \emph{Let $\calX$ be the desingularisation of the
  closure in $\PP^2_{\ZZ[\zeta_p]}$ of the Fermat curve $x^p+y^p=z^p$ 
with prime exponent $p$ as given in  \cite{Mc}.  
Then the arithmetic self-intersection
  number of its dualizing sheaf equipped with the Arakelov metric is
  bounded from above by
\begin{align*}
  \overline{\omega}_\Ar^2 &\le (2 g_p-2) \left( \log|\Delta_{\QQ(\zeta_p)|\QQ}|^2 + 
   [\QQ(\zeta_p):\QQ]  (  \kappa_1 \log p +\kappa_2)   
 +  g_p\, p^8 \log p     
  \right)
\end{align*}
where $\kappa_1, \kappa_2 \in \RR$ are absolute constants independent of $p$.}

{\bf Proof.}  Because of the previous discussion, we
are allowed to apply corollary \ref{cor:keyformula} with $P= \infty \in \PP^1$.
Then corollary \ref{cor:bound-bbelyi} and Theorem \ref{thm:geom-bound}
finish the proof. 
\hfill $\Box$
 
 \section{Bounds for  modular curves $X(N)$}
 
We now consider the modular curves
 $X(N)$, where $N=p^k m$ and we assume $p>3$ is a prime and $m\neq 1$ is coprime to $6 p$.
 These curves have a model over
 $\QQ(\zeta_N)$, where $\zeta_N$ is a primitive $N$-th root of unity,
 and its complex valued points correspond to the compact Riemann
 surfaces $\Gamma(N) \setminus \left( \HH \cup \PP^1(\QQ) \right)$,
 where $\Gamma(N)= \left\{ \left(\begin{smallmatrix} a & b \\ c &
       d\end{smallmatrix} \right) \in \Gamma(1) \,\left|
     \left(\begin{smallmatrix} a & b \\ c & d\end{smallmatrix} \right)
     \equiv \left(\begin{smallmatrix} 1 & 0 \\ 0 & 1\end{smallmatrix}
     \right) \mod{N} \right.  \right\}$. The natural morphism $X(N) \to X(1)$ 
is a Galois cover with three branch points. The ramification index above 
the cusp $\infty \in X(1)$ equals $N$.

 Regular models $\calX(N)$
over $\Spec \ZZ[\zeta_N]$ have been obtained by Deligne and Rapoport via normalisation 
and by Katz and Mazur via moduli interpretation. In the latter case
the scheme 
$\calX(N)$ is regular and it is smooth over $\ZZ[\zeta_N,1/N]$. If $\frakp|p$, then 
the fiber $\calX(N) \otimes \overline\FF_\frakp$
is the union of 
\begin{align*}
r_\frakp = p^k+p^{k-1}
\end{align*}
irreducible components crossing in the 
\begin{align*}
s_\frakp = \frac{p-1}{24} m^2 \phi(m) \prod_{q|m} \left(1+\frac{1}{q}\right)
\end{align*}
supersingular points of $\calX(N/p^k) \otimes
\overline\FF_\frakp$ with a maximal multiplicity  $m_\frakp = p^{2k-2}$.
From this description we deduce
\begin{align*}
u_\frakp \le (r_\frakp -1) m_\frakp s_\frakp , \quad l_\frakp \ge s_\frakp, \quad c_\frakp = 1
\end{align*}  
and therefore
\begin{align*}
 b_\frakp \le \frac{
 (r_\frakp - 1)^2 m_\frakp}{s_\frakp}.
\end{align*}
Finally we recall that the natural morphism $\belyi: \calX(N) \to \calX(1)$
is a Galois cover.

 \nn {\bf Theorem.} \emph{Let $\calX(N)$ be a  
 regular model of the modular curve $X(N)$ as above. 
Then the arithmetic self-intersection
  number of its dualizing sheaf equipped with the Arakelov metric is
  bounded from above by
\begin{align*}
  \overline{\omega}_\Ar^2 &\le (2 g_N-2) \left( \log|\Delta_{\QQ(\zeta_N)|\QQ}|^2 + 
   [\QQ(\zeta_N):\QQ]  (  \kappa_1 \log N +\kappa_2)   
 + 2  \sum_{\frakp \supset (N)}  \frac{(r_\frakp-1)^2 m_\frakp}{s_\frakp}  \log \Nm \frakp  
  \right),
\end{align*} 
where $\kappa_1, \kappa_2 \in \RR$ are absolute constants independent of $N$.}

{\bf Proof.}  
Because of the previous discussion, we
are allowed to apply corollary \ref{cor:keyformula} with $P= \infty \in \calX(1)$.
Now by means of corollary \ref{cor:bound-bbelyi} and Theorem \ref{thm:geom-bound}
the claim follows.
\hfill $\Box$

\end{document}